\documentclass{amsart}
\usepackage{amscd,amsmath,amssymb,amsfonts}
\usepackage[cmtip, all]{xy}

\newtheorem{thm}{Theorem}[section]
\newtheorem{prop}[thm]{Proposition}
\newtheorem{lem}[thm]{Lemma}
\newtheorem{cor}[thm]{Corollary}

\renewcommand{\theclaim}{\kern-3pt}

\newtheorem{IntroThm}{Theorem}

\theoremstyle{definition}
\newtheorem{Def}[thm]{Definition}

\theoremstyle{remark}
\newtheorem{rem}[thm]{Remark}
\newtheorem{rems}[thm]{Remarks}
\newtheorem{ex}[thm]{Example}

\newtheorem{IntroRem}{Remark}

\numberwithin{equation}{section}

\newcommand{\sC}{{\mathcal C}}

\newcommand{\sE}{{\mathcal E}}
\newcommand{\sF}{{\mathcal F}}
\newcommand{\sG}{{\mathcal G}}
\newcommand{\sH}{{\mathcal H}}
\newcommand{\sI}{{\mathcal I}}

\newcommand{\sN}{{\mathcal N}}
\newcommand{\sO}{{\mathcal O}}
\newcommand{\sP}{{\mathcal P}}

\newcommand{\sS}{{\mathcal S}}
\newcommand{\sT}{{\mathcal T}}
\newcommand{\sU}{{\mathcal U}}

\newcommand{\sX}{{\mathcal X}}
\newcommand{\sY}{{\mathcal Y}}

\newcommand{\A}{{\mathbb A}}

\newcommand{\G}{{\mathbb G}}

\renewcommand{\P}{{\mathbb P}}

\newcommand{\R}{{\mathbb R}}
\newcommand{\mS}{{\mathbb S}}

\newcommand{\Z}{{\mathbb Z}}

\renewcommand{\phi}{\varphi}

\newcommand{\codim}{{\rm codim}}

\newcommand{\Hom}{{\rm Hom}}
\newcommand{\End}{{\rm End}}

\newcommand{\Spec}{\operatorname{Spec}}

\newcommand{\Char}{\operatorname{char}}

\newcommand{\0}{\emptyset}
\newcommand{\sHom}{{\mathcal{H}{om}}}

\newcommand{\id}{{\operatorname{id}}}

\newcommand{\op}{{\text{\rm op}}}
\newcommand{\<}{\mathopen<}
\renewcommand{\>}{\mathclose>}
\newcommand{\Sets}{{\mathbf{Sets}}}

\renewcommand{\max}{{\operatorname{\rm max}}}

\newcommand{\Spt}{{\mathbf{Spt}}}
\newcommand{\Spc}{{\mathbf{Spc}}}
\newcommand{\Sm}{{\mathbf{Sm}}}

\newcommand{\hocolim}{\mathop{{\rm hocolim}}}
\renewcommand{\lim}{\operatornamewithlimits{\varprojlim}}
\newcommand{\colim}{\operatornamewithlimits{\varinjlim}}

\newcommand{\Tot}{{\operatorname{\rm Tot}}}

\newcommand{\SH}{{\operatorname{\sS\sH}}}

\newcommand{\eff}{{\mathop{eff}}}

\newcommand{\GW}{\operatorname{GW}}

\newcommand{\ds}{{/\kern-3pt/}}

\newcommand{\hofib}{{\mathop{\rm{hofib}}}}

\newcommand{\Bl}{\text{Bl}}

\newcommand{\Proj}{{\operatorname{Proj}}}

\newcommand{\Tate}{\text{Tate}}
\renewcommand{\:}{\kern-1.5pt:\kern-1.5pt}

\begin{document}

\title{The slice filtration and Grothendieck-Witt groups}
\author{Marc Levine}
\address{Universit\"at Duisburg-Essen\\
Fakult\"at Mathematik, Campus Essen\\
45117 Essen\\
Germany}
\email{marc.levine@uni-due.de}
\thanks{Research supported by the Alexander von Humboldt Foundation}


 
\renewcommand{\abstractname}{Abstract}
\begin{abstract}  Let $k$ be  a perfect field of  characteristic different from two. We show that the filtration on the Grothendieck-Witt group $\GW(k)$ induced by the slice filtration for  the sphere spectrum  in the motivic stable homotopy category  is the $I$-adic filtration, where $I$ is the augmentation ideal in $\GW(k)$. \end{abstract}
\date{\today}
\maketitle
\tableofcontents

\section*{Introduction}  Let $k$ be a perfect field of characteristic different from two. 
A fundamental theorem of Morel \cite{MorelA1, MorelLec} states that the endomorphism ring of the motivic sphere spectrum $\mS_k\in \SH(k)$ is naturally isomorphic to the Grothendieck-Witt ring of quadratic forms over $k$, $\GW(k)$. This result follows from Morel's calculation \cite[corollary 3.43]{MorelA1} of the corresponding bi-graded homotopy sheaves of $S^n\wedge \G_m^{\wedge q}$ in the unstable motivic homotopy category $\sH_\bullet(k)$ as the Milnor-Witt sheaves
\[
\pi_{m+p,p}(S^n\wedge \G_m^{\wedge {q}})\cong 
\begin{cases} \underline{K}^{MW}_{q-p}&\quad\text{ for } n=m\ge2, q\geq1, p\ge0,\\0&\quad\text{ for }m< n, p,q\ge0.\end{cases}
\]
Evaluating at $k$ and taking $m=n$, $p=q$ gives
\[
\End_{\sH_\bullet(k)}(S^{m}\wedge \G_m^{\wedge q})=K^{MW}_0(k)\text{ for }m\ge2, q\ge1.
\]
Combining this with Morel's isomorphism $K^{MW}_0(k)\cong \GW(k)$ and stabilizing gives Morel's theorem
\[
\End_{\SH(k)}(\mS_k)=\GW(k).
\]
This also leads to  the computation of the homotopy sheaf $\pi_{p,p}\Sigma^\infty_s\G_m^{\wedge {q}}$ (in the $S^1$-stable homotopy category $\SH_{S^1}(k)$) as $\underline{K}^{MW}_{q-p}$, for all $q\ge1$, $p\ge0$. 

In another direction, Voevodsky \cite{VoevSlice} has defined natural towers in $\SH(k)$ and $\SH_{S^1}(k)$, which are analogs of the classical Postnikov tower in $\SH$; we call each of these towers the {\em Tate Postnikov tower} (in $\SH(k)$ or  $\SH_{S^1}(k)$, as the case may be). Just as the classical Postnikov tower measures the $S^n$-connectivity of a spectrum, the Tate Postnikov tower measures the $S^{*,n}$ connectivity of a motivic spectrum.

In particular, the tower for $\mS_k$
\[
\ldots\to f_{n+1}\mS_k\to f_n\mS_k\to\ldots\to f_0\mS_k=\mS_k
\]
gives a filtration on the sheaf $\pi_{0,0}\mS_k$ by
\[
F^n_\Tate\pi_{0,0}\mS_k:=\text{im}(\pi_{0,0}f_n\mS_k\to \pi_{0,0}\mS_k).
\]
We have a similarly defined filtration on $\pi_{p,p}\Sigma_s^\infty\G_m^{\wedge q}$, which determines $F^n_\Tate\pi_{0,0}\mS_k:$ by
\[
F^n_\Tate\pi_{0,0}\mS_k:=\colim_qF_\Tate^{n+q}\pi_{q,q}\Sigma_s^\infty\G_m^{\wedge q}(k).
\]
 Our main result is the computation of $F^n_\Tate\pi_{p,p}\Sigma_s^\infty\G_m^{\wedge q}$, and thereby $F_\Tate^n\pi_{0,0}\mS_k$ (on perfect fields)
 
 \begin{IntroThm}\label{thm:main} Let $k$ be a perfect  field of  characteristic $\neq 2$ and let  $F$ be a perfect field extension of $k$. Let $n,p\ge0$, $q\ge1$ be integers and let
$N(a,b)=\max(0,\min(a,b))$.  Then via the identification given by Morel's isomorphism $\pi_{p,p}\Sigma_s^\infty\G_m^{\wedge q}\cong \underline{K}^{MW}_{q-p}$, we have
 \[
F^n_\Tate\pi_{p,p}\Sigma_s^\infty\G_m^{\wedge q}(F)= K^{MW}_{q-p}(F)\cdot I(F)^{N(n-p,n-q)},
 \]
 where $I(F)\subset  K^{MW}_{0}(F)$ is the  augmentation ideal. After stabilizing, this gives
 \[
F^n_\Tate\pi_{p,p}\Sigma^q_{\G_m}\mS_k(F) = K^{MW}_{q-p}(F)I(F)^{N(n-p,n-q)},\  n,p,q\in\Z,
 \]
 in particular,
 \[
F^n_\Tate\pi_{0,0}\mS_k(F) =  I(F)^{\max(n,0)}.
 \]
 \end{IntroThm}
See   theorem~\ref{thm:Main3}, corollary~\ref{cor:S1Stable} and corollary~\ref{cor:MainTStable}   for details. 

\begin{IntroRem} In case $k$ is a  field of characteristic $0$, we have a finer result, namely the identities stated in theorem~\ref{thm:main} extend to identities on the corresponding sheaves, for example
 \[
F^n_\Tate\pi_{p,p}\Sigma_s^\infty\G_m^{\wedge q} = \underline{K}^{MW}_{q-p} \cdot \sI^{N(n-p,n-q)}.
 \]
\end{IntroRem}

 Of course, one can more generally consider the filtration $F^*_\Tate\pi_{a,b}\sE$  on the homotopy sheaves $\pi_{a,b}\sE$ induced by the Tate Postnikov tower for an arbitrary $T$-spectrum $\sE\in\SH(k)$. In general, we cannot say anything about this filtration, but assuming a certain connectedness condition, we can compute the filtration on the first non-vanishing homotopy sheaves,  evaluated on perfect fields. 
  
 \begin{IntroThm}\label{thm:main0} Let $k$ be a perfect  field of   characteristic $\neq 2$ and let  $F$ be a perfect  field extension of $k$. Take $\sE\in \SH(k)$ and suppose that $\pi_{a+b,b}\sE=0$ for $a<0$, $b\in\Z$.  Then for $n> p$,
 \[
 F^n_\Tate\pi_{p,p}\sE(F) = [\pi_{n,n}\sE\cdot\underline{K}^{MW}_{n-p}]^{\widehat{\ }_{Tr}}(F).
\]
 For $n\le p$, we have  the identity of sheaves
\[
  F^n_\Tate\pi_{p,p}\sE=  \pi_{p,p}\sE.
\]
 \end{IntroThm}
 To explain the notation: The canonical action of $\pi_{*,*}\mS_k$ on $\pi_{*,*}\sE$, gives, for each finitely generated field extension $L$ of $k$, a right $K^{MW}_{-*}(L)$-module structure on $ \pi_{*,*}\sE(L)$, giving us the subgroup $\pi_{n,n}\sE(L)\cdot K^{MW}_{n-p}(L)$ of $\pi_{p,p}\sE(L)$.  This extends to arbitrary field extensions of $k$ by taking the evident colimit. Also, for each closed point $w\in \A^n_F$, we have a canonically defined {\em transfer map}
 \[
Tr_F(w)^*: \pi_{a,b}\sE(F(w))\to  \pi_{a,b}\sE(F)
\]
(see \S\ref{sec:Transfer} for details). $[\pi_{n,n}\sE\cdot\underline{K}^{MW}_{n-p}]^{\widehat{\ }_{Tr}}(F)$ is the subgroup of $\pi_{p,p}\sE(F)$ generated by the subgroups
$Tr_F(w)^*(\pi_{n,n}\sE(F(w))\cdot K^{MW}_{n-p}(F(w)))$, as $w$ runs over closed points of $\A_F^n$.  See  theorem~\ref{thm:Main2} for details.

Theorem~\ref{thm:main} is an easy consequence of theorem~\ref{thm:main0}; one uses Morel's unstable computations of the maps $S^{a,b}\wedge\Spec F_+\to S^{m,n}$ to reduce theorem~\ref{thm:main} to its $T$-stable version and then one uses the explicit presentation of $K^{MW}_*$  to compute
\[
[\underline{K}^{MW}_{q-n} \cdot \underline{K}^{MW}_{n-p}]]^{\widehat{\ }_{Tr}}(F)=K^{MW}_{q-p}(F)I^{N(n-p,n-q)}(F).
\]
 Morel's results on strictly $\A^1$-invariant sheaves allow us to go from the statement on functions fields to the one for the sheaves. 

The restriction to perfect fields arises from a separability assumption needed to compute the action of transfers on our selected generators for $F^n_\Tate\pi_{p,p}\sE$. We avoid characteristic two so as to have a description of the homotopy sheaves of the sphere spectrum in terms of Milnor-Witt $K$-theory. 

The paper is organized as follows. After setting up our notation and going over some background material on motivic homotopy theory in section \ref{sec:Background}, we recall some basic facts about the Tate Postnikov tower in section \ref{sec:HCT}. In section \ref{sec:Gen} we prove some connectedness results for the terms $f_nE$, $s_nE$ in the Tate Postnikov tower for an $S^1$-spectrum $E$ and give a description of generators for the subgroup $F^n_\Tate\pi_{0}E(F)$, all under a certain connectedness assumption on $E$. The generators are then factored into a product of two terms, one depending on $E$, the other only on the choice of a closed point of $\Delta^n_F\setminus\partial\Delta^n_F$. We analyze the second term in sections \ref{sec:PTCollapse}-\ref{sec:Transfer}, our main result being a description of this term as the $n$th suspension of a ``symbol map" associated to units $u_1,\ldots, u_n\in F^\times$. This is the main computation achieved in this paper. It is then relatively simple to feed this result into our description of the generators for  $F^n_\Tate\pi_{0}E(F)$ to prove  theorems \ref{thm:main}  and \ref{thm:main0} in section \ref{sec:Conclusion}; we conclude in section \ref{sec:Convergence}  with some remarks  on the convergence of the Tate Postnikov tower. 

I thank the referee for making several helpful suggestions and for pointing out a number of errors, including an incorrect formulation of theorem~\ref{thm:main0}, in an earlier version of this paper. Finally, I wish to thank the editors  for giving me the opportunity of contributing to this  volume. As a small token of my gratitude to Eckart for all of his aid and support over many years,   I dedicate this article to his memory.

\section{Background and notation} \label{sec:Background} Unless we specify otherwise, $k$ will be a fixed perfect base field, without restriction on the characteristic. For details on the following constructions, we refer the reader to \cite{GoerssJardine, Jardine, Jardine2, MorelA1, MorelConn,  MorelLec, MorelVoev}.

We write $[n]:=\{0,\ldots,n\}$ (including $[-1]=\0$) and let $\Delta$ be the category with objects $[n]$, $n=0,1,\ldots$, and morphisms $[n]\to[m]$ the order-preserving maps of sets. Given a category $\sC$, the category of simplicial objects in $\sC$ is as usual the category of functors $\Delta^\op\to\sC$.

$\Spc$ will denote the category of simplicial sets, $\Spc_\bullet$ the category of pointed simplicial sets, $\sH:=\Spc[WE^{-1}]$ the classical unstable homotopy category and $\sH_\bullet:=\Spc_\bullet[WE^{-1}]$ the pointed version. We denote the suspension operator $-\wedge S^1$ by $\Sigma_s$.  $\Spt$ is the category of suspension spectra and $\SH:=\Spt[WE^{-1}]$ the classical stable homotopy category.

The motivic versions are as follows: $\Sm/k$ is the category of smooth finite type $k$-schemes.  $\Spc(k)$ is the category of $\Spc$-valued presheaves on $\Sm/k$, $\Spc_\bullet(k)$ the $\Spc_\bullet$-valued presheaves, and $\Spt_{S^1}(k)$ the $\Spt$-valued presheaves. These all come with ``motivic"  model structures (see for example \cite{Jardine2}); we denote the corresponding homotopy categories by $\sH(k)$, $\sH_\bullet(k)$ and $\SH_{S^1}(k)$, respectively. Sending $X\in \Sm/k$ to the sheaf of sets  on $\Sm/k$ represented by $X$ (also denoted $X$)  gives an embedding of $\Sm/k$ to $\Spc(k)$; we have the similarly defined embedding of the category of smooth pointed schemes over $k$ into $\Spc_\bullet(k)$. All these categories are equipped with an internal Hom, denoted $\sHom$.

Let $\G_m$ be the pointed $k$-scheme $(\A^1\setminus0,1)$. In $\sH_\bullet(k)$ we have the objects $S^{a+b,b}:=\Sigma^a_s\G_m^{\wedge b}$, for $b\ge1$, $S^{n,0}:=S^n=\Sigma^n_s\Spec k_+$. If $X$ is a scheme with a $k$-point $x$, we write $(X,x)$ for the corresponding object in $\Spc_\bullet(k)$ or $\sH_\bullet(k)$. For a cofibration $\sY\to \sX$ in $\Spc(k)$, we usually give the quotient $\sX/\sY$ the canonical base-point $\sY/\sY$, but on occasion, we will give $\sX/\sY$ a base-point coming from a point $x\in \sX(k)$; we write this as $(\sX/\sY,x)$.

We let $T:=\A^1/(\A^1\setminus\{0\})$ and let $\Spt_T(k)$ denote the category of $T$-spectra, i.e.,  spectra in $\Spc_\bullet(k)$ with respect to the $T$-suspension functor $\Sigma_T:= -\wedge T$. $\Spt_T(k)$ has a motivic  model structure (see \cite{Jardine2}) and $\SH(k)$ is the homotopy category. We can also form the category of spectra in $\Spt_{S^1}(k)$ with respect to $\Sigma_T$;  with an appropriate model structure the resulting homotopy category is equivalent to $\SH(k)$. We will ignore the subtleties of this distinction and simply identify the two homotopy categories.

Both $\SH_{S^1}(k)$ and $\SH(k)$ are triangulated categories with suspension functor $\Sigma_s$. We have the   triangle of {\em infinite suspension functors} $\Sigma^\infty$ and their right adjoints $\Omega^\infty$
\[
\xymatrix{
\sH_\bullet(k)\ar[r]^{\Sigma^\infty_s}\ar[rd]_{\Sigma^\infty_T}&\SH_{S^1}(k)\ar[d]^{\Sigma^\infty_T}\\
&\SH(k)
}\quad
\xymatrix{
\sH_\bullet(k)&\SH_{S^1}(k)\ar[l]_{\Omega^\infty_s}\\
&\SH(k)\ar[u]^{\Omega^\infty_T}\ar[ul]^{\Omega^\infty_T}
}
\]
both commutative up to natural isomorphism. These are all left, resp. right derived versions of Quillen adjoint pairs of functors on the underlying model categories. We note that the suspension functor $\Sigma_{\G_m}$ is invertible on $\SH(k)$.

For $\sX\in \sH_\bullet(k)$, we have the bi-graded homotopy sheaf $\pi_{a,b}\sX$, defined for $b\ge0$, $a-b\ge0$, as the Nisnevich sheaf associated to the presheaf on $\Sm/k$
\[
U\mapsto \Hom_{\sH_\bullet(k)}(\Sigma^{a-b}_s\Sigma^b_{\G_m}U_+,\sX).
\]
These extend in the usual way to bi-graded homotopy sheaves $\pi_{a,b}E$ for $E\in\SH_{S^1}(k)$, $b\ge0$, $a\in\Z$, and $\pi_{a,b}\sE$ for $\sE\in\SH(k)$,  $a,b\in\Z$, by taking the Nisnevich sheaf associated to 
\[
U\mapsto \Hom_{\SH_{S^1}(k)}(\Sigma^{a-b}_s\Sigma^b_{\G_m}\Sigma^\infty_sU_+,E)
\text{ or }
U\mapsto \Hom_{\SH(k)}(\Sigma^{a-b}_s\Sigma^b_{\G_m}\Sigma^\infty_TU_+,\sE),
\]
as the case may be. We write $\pi_n$ for $\pi_{n,0}$; for e.g. $E\in \Spt_{S^1}(k)$ fibrant, $\pi_nE$ is the Nisnevich sheaf associated to the presheaf $U\mapsto \pi_n(E(U))$.

For $F$ a finitely generated field extension of $k$, we may view $\Spec F$ as the generic point of some $X\in\Sm/k$. Thus, for a Nisnevich sheaf $\sS$ on $\Sm/k$, we may define $\sS(F)$ as the stalk of $\sS$ at $\Spec F\in X$. For an arbitrary field extension $F$  of $k$ (not necessarily finitely generated over $k$), we define $\sS(F)$  as the colimit over $\sS(F_\alpha)$, as $F_\alpha$ runs over subfields of $F$ containing $k$ and finitely generated over $k$.

 \section{The homotopy coniveau tower}\label{sec:HCT}  Our computations rely heavily on our model for the Tate Postnikov tower in $\SH_{S^1}(k)$, which we briefly recall (for details, we refer the reader to \cite{LevineHC}). We start by recalling the Tate Postnikov tower in $\SH_{S^1}(k)$ and introducing some notation.
 
Fix a perfect base-field $k$. Let 
\[
\Sigma_T:\SH_{S^1}(k)\to \SH_{S^1}(k)
\]
be the $T$-suspension functor. For $n\ge0$, we let $\Sigma_T^n\SH_{S^1}(k)$ be the localizing subcategory of $\SH_{S^1}(k)$ generated by infinite suspension spectra of the form $\Sigma^n_T\Sigma^\infty_sX_+$, with $X\in \Sm/k$. We note that $\Sigma^0_T\SH_{S^1}(k)=\SH_{S^1}(k)$. The inclusion functor $i_n:\Sigma_T^n\SH_{S^1}(k)\to \SH_{S^1}(k)$ admits, by results of Neeman \cite{Neeman}, a right adjoint $r_n$; define the functor
$f_n:\SH_{S^1}(k)\to\SH_{S^1}(k)$ by $f_n:=i_n\circ r_n$. The unit for the adjunction gives us the natural morphism
\[
\rho_n:f_nE\to E
\]
for $E\in \SH_{S^1}(k)$; similarly, the inclusion $\Sigma_T^m\SH_{S^1}(k)\subset \Sigma_T^n\SH_{S^1}(k)$ for $n<m$ gives the natural transformation  $f_mE\to f_nE$, forming the {\em Tate Postnikov tower}
\[
\ldots\to f_{n+1}E\to f_nE\to\ldots\to f_0E=E.
\]
We complete $f_{n+1}E\to f_nE$ to a distinguished triangle
\[
f_{n+1}E\to f_nE\to s_nE\to f_{n+1}E[1];
\]
this turns out to be functorial in $E$. The object $s_nE$ is the {\em $n$th slice} of $E$.

There is an analogous construction in $\SH(k)$: For $n\in\Z$, let $\Sigma^n_T\SH^\eff(k)\subset\SH(k)$ be the localizing category generated by the $T$-suspension spectra $\Sigma^n_T\Sigma^\infty_TX_+$, for $X\in\Sm/k$. As above, the inclusion $i_n:\Sigma^n_T\SH^\eff(k)\to\SH(k)$ admits a left adjoint $r_n$, giving us the truncation functor $f_n$ and the Postnikov tower
\[
\ldots\to f_{n+1}\sE\to f_n\sE\to\ldots\to \sE.
\]
Note that this tower is in general infinite in both directions. We define the layer $s_n\sE$ as above.

By \cite[theorem 7.4.1]{LevineHC}, the 0-space functor $\Omega^\infty_T$ sends $\Sigma^n_T\SH^\eff(k)$ to $\Sigma_T^n\SH_{S^1}(k)$. This fact, together with the universal properties of the truncation functors $f_n$ in $\SH_{S^1}(k)$ and $\SH(k)$, plus the fact that $\Omega^\infty_T$ is a right adjoint,  gives the canonical isomorphism for $n\ge0$
\begin{equation}\label{eqn:SliceIso}
f_n\Omega^\infty_T\sE\cong \Omega^\infty_Tf_n\sE.
\end{equation}
Furthermore, for $E\in \SH_{S^1}(k)$, we have (by  \cite[theorem 7.4.2]{LevineHC}) the canonical isomorphism
\begin{equation}\label{eqn:SliceIso2}
\Omega_{\G_m}f_nE=f_{n-1}\Omega_{\G_m}E.
\end{equation}
As $\Omega_{\G_m}:\SH(k)\to \SH(k)$ is an auto-equivalence,  and restricts to an equivalence
\[
\Omega_{\G_m}:\Sigma^n_T\SH^\eff(k)\to \Sigma^{n-1}_T\SH^\eff(k),
\]
the analogous identity in $\SH(k)$ holds as well.

\begin{Def} For $a\in \Z$, $b\ge0$, $E\in \SH_{S^1}(k)$, define the filtration $F^n_\Tate\pi_{a,b}E$, $n\ge0$, of $\pi_{a,b}E$  by
\[
F^n_\Tate\pi_{a,b}E:=\text{im}(\pi_{a,b}f_nE\to  \pi_{a,b}E).
\]
Similarly, for $\sE\in \SH(k)$, $a,b,n\in\Z$, define
\[
F^n_\Tate\pi_{a,b}\sE:=\text{im}(\pi_{a,b}f_n\sE\to  \pi_{a,b}\sE).
\]
\end{Def}
The main object of this paper is to understand $F^n_\Tate\pi_0E$ for suitable $E$. For later use, we note the following

\begin{lem} \label{lem:DegreeShift} 1. For $E\in\SH_{S^1}(k)$, $n,p,a,b\in\Z$ with $n,p, b, n-p, b-p\ge0$,  the adjunction isomorphism $\pi_{a,b}E\cong \pi_{a-p.b-p}\Omega^p_{\G_m}E$ induces an isomorphism
\[
F^n_\Tate\pi_{a,b}E\cong F^{n-p}_\Tate\pi_{a-p,b-p}\Omega^p_{\G_m}E.
\]
Similarly, for $\sE\in \SH(k)$, $n,p,a,b\in \Z$,  the adjunction isomorphism $\pi_{a,b}\sE\cong \pi_{a-p.b-p}\Omega^p_{\G_m}\sE$ induces an isomorphism
\[
F^n_\Tate\pi_{a,b}\sE\cong F^{n-p}_\Tate\pi_{a-p,b-p}\Omega^p_{\G_m}\sE.
\]
2. For $\sE\in \SH(k)$, $a,b,n\in\Z$, with $b,n\ge0$, we have a canonical isomorphism
\[
\phi_{\sE,a,b,n}: \pi_{a,b}f_n\sE\to \pi_{a,b}\Omega^\infty_Tf_n\sE,
\]
inducing an isomorphism $F^n_\Tate\pi_{a,b}\sE\cong F^n_\Tate\pi_{a,b}\Omega_T^\infty\sE$.
\end{lem}

\begin{proof} (1) By \eqref{eqn:SliceIso2}, adjunction induces isomorphisms
\begin{align*}
F^n_\Tate\pi_{a,b}E&:=\text{im}(\pi_{a,b}f_nE\to \pi_{a,b}E)\\
&\cong \text{im}(\pi_{a-p,b-p}\Omega^p_{\G_m}f_nE\to \pi_{a-p,b-p}\Omega^p_{\G_m}E)\\
&=  \text{im}(\pi_{a-p,b-p} f_{n-p}\Omega^p_{\G_m}E\to \pi_{a-p,b-p}\Omega^p_{\G_m}E)\\
&=F^{n-p}_\Tate\pi_{a-p,b-p}\Omega^p_{\G_m}E.
\end{align*}
The proof for $\sE\in \SH(k)$ is the same.

For (2), the isomorphism $\phi_{\sE,a,b,n}$ arises from \eqref{eqn:SliceIso} and the adjunction isomorphism
\begin{align*}
\Hom_{\SH_{S^1}(k)}(\Sigma^\infty_s\Sigma^{a-b}_s\Sigma^b_{\G_m}U_+, f_n\Omega^\infty_T\sE)&\cong
\Hom_{\SH_{S^1}(k)}(\Sigma^\infty_s\Sigma^{a-b}_s\Sigma^b_{\G_m}U_+, \Omega^\infty_Tf_n\sE)\\
&\cong
\Hom_{\SH(k)}(\Sigma^\infty_T\Sigma^{a-b}_s\Sigma^b_{\G_m}U_+, \sE).
\end{align*}
\end{proof}

We now turn to a discussion of our model for $f_nE(X)$, $X\in\Sm/k$. We start with the cosimplicial scheme $n\mapsto \Delta^n$, with $\Delta^n$ the {\em algebraic $n$-simplex} $\Spec k[t_0,\ldots, t_n]/\sum_it_i-1$. The cosimplicial structure is given by sending a map $g:[n]\to [m]$ to the map $g:\Delta^n\to \Delta^m$ determined by
 \[
 g^*(t_i)=\begin{cases}\sum_{j, g(j)=i}t_j&\text{ if }g^{-1}(i)\neq\0\\0&\text{ else.}\end{cases}
 \]
 A {\em face} of $\Delta^m$ is a closed subscheme $F$ defined by equations $t_{i_1}=\ldots=t_{i_r}=0$; we let $\partial\Delta^n\subset \Delta^n$ be the closed subscheme defined by $\prod_{i=0}^nt_i=0$, i.e.,  $\partial\Delta^n$ is the union of all the proper faces.

 Take $X\in \Sm/k$. We let $\sS_X^{(q)}(m)$ denote the set of closed subsets $W\subset X\times\Delta^m$ such that $\codim_{X\times F}W\cap X\times F\ge q$ for all faces $F\subset \Delta^m$ (including $F=\Delta^m$). We make $\sS_X^{(q)}(m)$ into a partially ordered set via inclusions of closed subsets. Sending $m$ to  $\sS_X^{(q)}(m)$ and $g:[n]\to [m]$ to $g^{-1}:\sS_X^{(q)}(m)\to \sS_X^{(q)}(n)$ gives us the simplicial poset $\sS_X^{(q)}$. 
 
 Now take $E\in \Spt_{S^1}(k)$. For $X\in\Sm/k$ and closed subset $W\subset X$, we have the spectrum with supports $E^W(X)$ defined as the homotopy fiber of the restriction map $E(X)\to E(X\setminus W)$. This construction is functorial in the pair $(X,W)$, where we define a map $f:(Y,T)\to (X,W)$ as a morphism $f:Y\to X$ in $\Sm/k$ with $f^{-1}(W)\subset T$. 
 
 Define 
 \[
 E^{(q)}(X,m):=\hocolim_{W\in \sS_X^{(q)}(m)}E^W(X\times\Delta^m).
 \]
 The fact that $m\mapsto \sS_X^{(q)}(m)$ is a simplicial poset, and $(Y,T)\mapsto E^T(Y)$ is a functor from the category of pairs to spectra shows that $m\mapsto  E^{(q)}(X,m)$ defines a simplicial spectrum. We denote the associated total spectrum by $E^{(q)}(X)$.

 For $q\ge q'$, the inclusions $\sS_X^{(q)}(m)\subset \sS_X^{(q')}(m)$ induces a map of simplicial posets $\sS_X^{(q)}\subset \sS_X^{(q')}$ and thus a morphism of spectra $i_{q',q}:E^{(q)}(X)\to E^{(q')}(X)$. We have as well the natural map 
 \[
 \epsilon_X:E(X)\to \Tot(E(X\times\Delta^*))= E^{(0)}(X), 
 \]
 which is a weak equivalence if $E$ is homotopy invariant. Together, this forms the {\em augmented homotopy coniveau tower}  tower
 \[
E^{(*)}(X):=\ldots\to E^{(q+1)}(X)\xrightarrow{i_{q}}E^{(q)}(X)\xrightarrow{i_{q-1}}\ldots E^{(1)}(X)
 \xrightarrow{i_{0}}E^{(0)}(X)\xleftarrow{\epsilon_X}E(X)
\]
with $i_q:=i_{q,q+1}$. Thus, for homotopy invariant $E$, we have the 
homotopy coniveau tower in $\SH$
 \[
E^{(*)}(X):= \ldots\to E^{(q+1)}(X)\xrightarrow{i_{q}}E^{(q)}(X)\xrightarrow{i_{q-1}}\ldots E^{(1)}(X)
 \xrightarrow{i_{0}}E^{(0)}(X)\cong E(X).
\]

Letting $\Sm\ds k$ denote the subcategory of $\Sm/k$ with the same objects and with morphisms the smooth morphisms, it is not hard to see that sending $X$ to $E^{(*)}(X)$ defines a functor from  $\Sm\ds k^\op$ to augmented towers of spectra.

 On the other hand, for $E\in \Spt_{S^1}(k)$, we have the (augmented) Tate Postnikov tower
 \[
 f_*E:=\ldots\to f_{q+1}E\to f_qE\to\ldots\to f_0E\cong E
 \]
 in $\SH_{S^1}(k)$, which we may evaluate at $X\in\Sm/k$, giving the  tower $f_*E(X)$ in  $\SH$, augmented over $E(X)$.
 
As a direct consequence of  our main result (theorem 7.1.1) from \cite{LevineHC} we have
 \begin{thm} \label{thm:HC} Let $E$ be a quasi-fibrant object in $\Spt_{S^1}(k)$ for the model structure described in \cite{Jardine}, and take $X\in \Sm/k$. Then there is an isomorphism of augmented towers in $\SH$
 \[
( f_*E)(X)\cong E^{(*)}(X)
\]
over the identity on $E(X)$, which is natural with respect to smooth morphisms in $\Sm/k$.
 \end{thm}
 
  In particular, we may use the explicit model $E^{(q)}(X)$ to understand $(f_qE)(X)$.

\begin{rem} For $X, Y\in \Sm/k$ with given $k$-points $x\in X(k)$, $y\in Y(k)$, we have a natural isomorphism in $\SH_{S^1}(k)$
\[
\Sigma^\infty_s(X\wedge Y)\oplus \Sigma^\infty_s(X\vee Y)\cong \Sigma^\infty_s(X\times Y)
\]
i.e. $\Sigma^\infty_s(X\wedge Y)$ is a canonically defined summand of $\Sigma^\infty_s(X\times Y)$. In particular for $E$ a quasi-fibrant object of $\Spt_{S^1}(k)$, we have a natural isomorphism in $\SH$
\[
\sHom(\Sigma^\infty_s(X\wedge Y),E)\cong \hofib\left(E(X\times Y)\to \hofib(E(X)\oplus E(Y)\to E(k))\right)
\]
where the maps are induced by the evident restriction maps. In particular, we may define $E(X\wedge Y)$ via the above isomorphism, and our comparison results for Tate Postnikov tower and homotopy coniveau tower extend to values at smash products of smooth pointed schemes over $k$.
\end{rem}

\section{Connectedness and generators for $\pi_0$}\label{sec:Gen} As in section \ref{sec:HCT}, our base-field $k$ is perfect.  We fix a quasi-fibrant $S^1$-spectrum $E\in\Spt_{S^1}(k)$.

\begin{lem}\label{lem:Vanishing} Let $F$ be a finitely generated  field extension of $k$, $x\in\A^n_F$ a closed point. Then for every $m>0$, the map
\[
i_{0*}:E^{(x,0)}(\A^n\times\A^m_F)\to E^{(x\times_F \A^m_F)}(\A^n\times\A^m_F)
\]
induced by the map of pairs 
\[
\id_{\A^n\times\A^m}:(\A^n\times\A^m_F, x\times\A^m_F)\to
(\A^n\times\A^m_F, (x,0))
\]
is the zero-map in $\SH$. In particular, the induced map on homotopy groups is the zero map.
\end{lem}

\begin{proof} We use the Morel-Voevodsky purity isomorphisms in $\sH_\bullet(k)$  \cite[Theorem 3.2.23]{MorelVoev}, with the isomorphisms defined via a fixed choice of generators for the maximal ideal $m_x\subset \sO_{\A^n_F,x}$ and $m_0\subset \sO_{\A^m,0}$
\begin{align*}
\A^n_F\times\A^m/(\A^n_F\times\A^m\setminus\{(x,0)\})&\cong\Sigma_T^{n+m} (x,0)_+\\
&\cong \Sigma^n_Tx\times\A^m/(x\times\A^m\setminus\{(x,0)\})\\
\A^n_F\times\A^m/(\A^n_F\times\A^m\setminus x\times\A^m)&\cong\Sigma_T^{n} x\times \A^m_{+}.
\end{align*}
Via these isomorphisms, the quotient map 
\[
q:\A^n_F\times\A^m/(\A^n_F\times\A^m\setminus x\times\A^m)\to \A^n_F\times\A^m/(\A^n_F\times\A^m\setminus \{(x,0\})
\]
is isomorphic to the $n$th $T$-suspension of the quotient map
\[
q':x\times\A^m_{+}\to x\times\A^m/(x\times\A^m\setminus\{(x,0)\})
\]
As $i_{0*}$ is the map  induced by applying $\sHom(-,E)$ to $\Sigma^\infty_s q$, we need only show that $q'$ factors through the map $x\times\A^m_{+}\to *$ (in $\sH_\bullet(k)$). This follows from the commutative diagram
\[
\xymatrix{
x\times 1_+\ar[r]^i_\sim\ar[d]&x\times\A^m_{+}\ar[d]^{q'}\\
{*}\ar[r]&x\times\A^m/(x\times\A^m\setminus\{(x,0)\}),
}
\]
where $1=(1,\ldots, 1)\in \A^m$, since $i$ is an isomorphism in $\sH_\bullet(k)$ by homotopy invariance.
\end{proof}

We have the re-indexed homotopy sheaves  
$\Pi_{n,m}(E):=\pi_{n+m,m}(E)$.
We have as well the sheaf $\pi_nE:=\pi_{n,0}E$; we call $E$ $m$-connected if $\pi_n(E)=0$ for all $n\le m$.

 Since $E^{(n)}(X)=\Tot[m\mapsto E^{(n)}(X,m)]$, we have the  strongly convergent spectral sequence
\begin{equation}\label{eqn:SimpSSeq}
E^1_{p,q}(X)=\pi_{q}E^{(n)}(X,p)\Longrightarrow \pi_{p+q}E^{(n)}(X),
\end{equation}
Now take $X=\Spec F$, $F$ a finitely generated  field over $k$. For dimensional reasons, we have $\sS^{(n)}_F(p)=\0$ for $p<n$, and we therefore have an edge homomorphism
\[
\epsilon_{-n}:\pi_{q-n}E^{(n)}(X,n)\to\pi_{q}E^{(n)}(X).
\]
Furthermore,  $\sS^{(n)}_F(n)$ is the set of closed points $w\in\Delta_F^n\setminus\partial\Delta^n_F$, so $\epsilon_{-n}$ can be written as
\[
\epsilon_{-n}:\oplus_{w\in(\Delta_F^n\setminus\partial\Delta^n_F)^{(n)}}\pi_{q-n}E^w(\Delta^n_F)\to \pi_qE^{(n)}(F);
\]
here $Y^{(a)}$ denotes the set of codimension $a$ points on a scheme $Y$.

Via the weak equivalence $E^{(n)}(F)\cong f_nE(F)$, we have the canonical map
\[
\epsilon_{-n}:\oplus_{w\in(\Delta_F^n\setminus\partial\Delta^n_F)^{(n)}}\pi_{q-n}E^w(\Delta^n_F)\to \pi_{q}f_nE(F).
\]
Similarly, composing with $f_nE\to s_nE$, we have the canonical map
\[
\epsilon_{-n}:\oplus_{w\in(\Delta_F^n\setminus\partial\Delta^n_F)^{(n)}}\pi_{q-n}E^w(\Delta^n_F)\to \pi_qs_nE(F).
\]

\begin{prop}\label{prop:Connected} Let $E\in\Spt_{S^1}(k)$ be quasi-fibrant. Suppose $\Pi_{a,*}E(F)=0$ for all $a<0$ and for all finitely generated  field extensions $F$ of $k$. Then for $n\ge0$: \\\\
1. $\Pi_{a,*}f_nE$ and $\Pi_{a,*}s_nE$ are zero for all $a<0$. In particular,  $f_nE$ and $s_nE$ are $-1$-connected. \\
\\
2. For each finitely generated  field $F$ over $k$, the edge homomorphisms
\[
\epsilon_{-n}:\oplus_{w\in(\Delta_F^n\setminus\partial\Delta^n_F)^{(n)}}\pi_{-n}E^w(\Delta^n_F)\to \pi_0(f_nE)(F)
\]
\[
\epsilon_{-n}:\oplus_{w\in(\Delta_F^n\setminus\partial\Delta^n_F)^{(n)}}\pi_{-n}E^w(\Delta^n_F)\to \pi_0(s_nE)(F)
\]
are surjections.
\end{prop}

\begin{proof}  Using the distinguished triangle
\[
f_{n+1}E\to f_nE\to s_nE\to \Sigma_s f_{n+1}E
\]
we see that it suffices to prove the statements for $f_nE$. 

Using the isomorphism \eqref{eqn:SliceIso2}, we see that for (1),  it suffices to show that $f_nE$ is $-1$-connected. By a theorem of Morel  \cite[lemma 3.3.6]{MorelLec}, it suffices to show that $f_nE(F)$ is $-1$-connected for all finitely generated field extensions $F$ of $k$. 

We first show that, for each $p\ge n$,
\begin{align*}
\text{a. } &\pi_{q}E^{(n)}(F,p)=0\text{ for } q<-p\\
\text{b. } &\text{The natural map }\\
&\hskip20pt\oplus_{W\in \sS^{(n)}_F(p), w\in W\cap (\Delta^p_F)^{(p)}}\pi_{-p}E^w(\Delta^p_F)\to 
\pi_{-p}E^{(n)}(F,p)\\
&\text{ is surjective.}\\
\end{align*}

For (a), let $W\subset \Delta^p_F$ be a closed subset. We have the Gersten spectral sequence
\[
E_1^{a,b}=\oplus_{w\in W\cap (\Delta^p_F)^{(a)}}\pi_{-a-b}E^w(\Spec \sO_{\Delta^p_F,w})\Longrightarrow
\pi_{-a-b}E^W(\Delta^p_F).
\]
Since $E$ is quasi-fibrant, and $\Delta^p_F$ is smooth over $k$, we have an isomorphism  (via Morel-Voevodsky purity \cite[Theorem 3.2.23]{MorelVoev})
\[
\pi_m(E^w(\Spec \sO_{\Delta^p_F,w}))\cong \pi_m(E(w_+\wedge S^{2a,a})),
\]
where $a=\codim_{\Delta^p_F}w$. But 
\[
\pi_m(E(w_+\wedge S^{2a,a}))=(\pi_{m+2a,a}E)(F(w))
\]
which is zero for $m+a<0$. Since $0\le a\le p$, we see that, for $m<-p$,
\[
\pi_mE^W(\Delta^p_F)=0.
\]
 As $E^{(n)}(F,p)$ is a colimit over $E^W(\Delta^p_F)$ with $W\in \sS^{(n)}_F(p)$, it follows that $\pi_mE^{(n)}(F,p)=0$ for $m<-p$, proving (a).

The same computation shows that $\pi_{-p}(E^w(\Spec \sO_{\Delta^p_F,w}))=0$ if $\codim_{\Delta^p_F}w<p$, so (b) follows from the Gersten spectral sequence.

Using the strongly convergent spectral sequence \eqref{eqn:SimpSSeq}, 
we see that (a) implies that $\pi_qE^{(n)}(F)=0$ for $q<0$.

Next, we show that
\begin{align*}
\text{c. }&\pi_{-p}E^{(n)}(F,p)=0\text{ for }p>n.
\end{align*}
For this, it suffices by (b) to show that for $w\in W\cap (\Delta_F^p)^{(p)}$ with $W\in \sS^{(n)}_F(p)$ and with  $p>n$, the map
\begin{equation}\label{eqn:ZeroMap}
\pi_{-p}E^w(\Delta^p_F)\to \pi_{-p}E^{(n)}(F,p)
\end{equation}
is the zero map. To see this, note that $W$ does not intersect any face $T$ of $\Delta^p_F$ having $\dim_FT< n$. Thus, there is a linear $W'\cong\A^{p-n}_{F'}\subset \Delta^p_F$ containing $w$ (for $F'$ some extension field of $F$ contained in $F(w)$) with $W'\in \sS^{(n)}_F(p)$: for a suitable degeneracy map $\sigma:\Delta^p\to \Delta^n$ one takes $W'=\sigma^{-1}(\sigma(w))$. By lemma~\ref{lem:Vanishing}, the map $E^w(\Delta^p_F)\to E^{W'}(\Delta^p_F)$ is the zero map in $\SH$; passing to the limit over all $W''\in \sS_F^{(n)}(p)$, we see that \eqref{eqn:ZeroMap} is the zero map, as claimed.

In the spectral sequence \eqref{eqn:SimpSSeq}, we have $E_{p,-p}^1=0$ for $p>n$; we also have 
$E_{p,-p}^1=0$ for $p<n$ since $\sS^{(n)}_F(p)=\0$ if $p<n$ for dimensional reasons. Thus, the only term contributing to $\pi_0E^{(n)}(F)$ is $E^1_{n,-n}$.   As the spectral sequence is strongly convergent, the edge homomorphism in the spectral sequence \eqref{eqn:SimpSSeq} induces a surjection
\[
\oplus_{w\in \sS^{(n)}_F(n)}\pi_{-n}E^w(\Delta^n_F)\to\pi_0E^{(n)}(F).
\]
Combining this with theorem~\ref{thm:HC} gives us the surjection 
\[
\oplus_{w\in \sS^{(n)}_F(n)}\pi_{-n}E^w(\Delta^n_F)\to\pi_0(f_nE(F)).
\]
Similarly, the vanishing $\pi_pE^{(n)}(F)=0$ for $p<0$  shows that $f_nE(F)$ is -1 connected. 
\end{proof}

We thus have generators $\oplus_{w\in(\Delta_F^n\setminus\partial\Delta^n_F)^{(n)}}\pi_{-n}E^w(\Delta^n_F)$ for $\pi_0f_nE(F)$, and hence for our main object of study, $F^n_\Tate\pi_0E(F)$. We examine the composition
\begin{equation}\label{eqn:*}
\pi_{-n}E^w(\Delta^n_F) \xrightarrow{\epsilon_{-n}} \pi_0f_nE(F)\xrightarrow{\rho_n}\pi_0E(F)
\end{equation}
more closely.

Fix  a closed point $w$ in $\Delta_F^n\setminus\partial\Delta^n_F$. We have the quotient map
\[
c_w:\Delta^n_F/\partial\Delta^n_F\to  \Delta^n_F/(\Delta^n_F\setminus w)
\]
and the canonical identification
\[
E^w(\Delta^n_F)=\sHom(\Sigma^\infty_s\Delta^n_F/(\Delta^n_F\setminus w),E).
\]
Thus, given an element $\tau\in \pi_{-n}(E^w(\Delta^n_F))$, we have the corresponding morphism
\[
\tau:\Sigma^\infty_s\Delta^n_F/(\Delta^n_F\setminus w)\to\Sigma_s^nE
\]
and we may compose with $c_w$ to give the map
\[
\tau\circ \Sigma^\infty_s c_w:\Sigma^\infty_s\Delta^n_F/\partial\Delta^n_F\to \Sigma_s^nE.
\]

As each of the faces of $\Delta^n_F$ are affine spaces over $F$, we have a canonical isomorphism
\[
\sigma_F:\Sigma^n_s\Spec F_+\to \Delta^n_F/\partial\Delta^n_F
\]
in $\sH_\bullet(k)$ (see the beginning of  \S\ref{sec:PTCollapse} for details), giving us the element
\[
\pi(\tau):=\tau\circ \Sigma^\infty_s(c_w\circ\sigma_F)\in \pi_n( \Sigma_s^nE(F))=\pi_0(E(F)).
\]
 
The following result is a direct consequence of the definitions:

\begin{lem} For $\tau\in \pi_{-n}(E^w(\Delta^n_F))$, $\pi(\tau)=\rho_n(\epsilon_{-n}(\tau))$. 
\end{lem}

On the other hand, we have the  Morel-Voevodsky purity isomorphism ({\it loc. cit.})
\begin{equation}\label{eqn:MV}
MV_w: \Delta^n_F/(\Delta^n_F\setminus w)\to w_+\wedge S^{2n,n}.
\end{equation}
The definition of $MV_w$ requires some additional choices;  we complete our definition of $MV_w$  in \S\ref{sec:Susp1}, where it is written as
$MV_w=(\id_{w_+}\wedge \alpha)\circ mv_w$ (see definition~\ref{Def:PurityIsomv_w} and \eqref{eqn:MVDef}).

In any case,  via $MV_w$,  we may factor $\pi(\tau)$ as 
\begin{align*}
\pi(\tau)&:= \tau\circ \Sigma^\infty_s(c_w\circ\sigma_F)\\
&= (\tau\circ  \Sigma^\infty_s MV_w^{-1})\circ   \Sigma^\infty_s(MV_w\circ c_w\circ\sigma_F)
\end{align*}
The term $\tau\circ  \Sigma^\infty_s MV_w^{-1}$ is the morphism
\[
\tau\circ  \Sigma^\infty_s MV_w^{-1}: \Sigma^\infty_s w_+\wedge S^{2n,n}\to \Sigma_s^nE
\]
which we may interpret as an element of $\pi_{-n}(\Omega_T^nE(w))$, while the morphism $\Sigma^\infty_s(MV_w\circ c_w\circ\sigma_F)$ is the infinite suspension of the map
\begin{equation}\label{eqn:CollapseMap}
Q_F(w):=MV_w\circ c_w\circ\sigma_F:\Sigma^n_s\Spec F_+\to w_+\wedge S^{2n,n}.
\end{equation}
Conversely, given any element 
$\xi\in \pi_{-n}(\Omega_T^nE(w))$, which we write as a morphism
\[
\xi:w_+\wedge S^{2n,n}\to \Sigma^n_sE
\]
we recover an element $\tau\in  \pi_{-n}(E^w(\Delta^n_F))$ as $\tau:=\xi\circ  \Sigma^\infty_s MV_w$, and thus the element
\[
\xi\circ\Sigma^\infty_sQ_F(w)\in \pi_n(\Sigma^n_sE(F))=\pi_0E(F)
\]
is in $F^n_\Tate\pi_0E(F)$.

Putting this all together, we have
\begin{prop}\label{prop:Factorization} Let $F$ be a finitely generated  field extension of $k$ and let $E\in \Spt_{S^1}(k)$ be quasi-fibrant.\\
\\
1. Let $w$ be a closed point of $\Delta^n_F\setminus\partial\Delta^n_F$, and take $\xi_w\in \pi_{-n}(\Omega_T^nE(w))$. Then 
$\xi_w\circ\Sigma^\infty_sQ_F(w)$ is in $F^n_\Tate\pi_0E(F)$.\\
\\
2. Suppose that $\Pi_{a,*}E=0$ for all $a<0$. Then  $F^n_\Tate\pi_0E(F)$ is generated by elements of the form $\xi_w\circ\Sigma^\infty_sQ_F(w)$, $\xi_w\in \pi_{-n}(\Omega_T^nE(w))$, as $w$ runs over closed points of $\Delta^n_F\setminus\partial\Delta^n_F$.
\end{prop}

\begin{rem} The proposition extends without change to arbitrary field extensions $F$ of $k$, by a simple limit argument.
\end{rem}

The next few sections will be devoted to giving explicit formulas for the map $Q_F(w)$. In case $w$ is an $F$-point of $\Delta^n\setminus\partial \Delta^n$, we are able to do so directly; in general, we will need to pass to an $n$-fold $\P^1$-suspension before we can give an explicit formula. We will then conclude with the proof of our main result in \S\ref{sec:Conclusion}.

\section{The Pontryagin-Thom collapse map}\label{sec:PTCollapse} We recall a special case of Pontryagin-Thom construction in $\sH_\bullet(k)$.

Let  $V_n$ be the open subscheme $\Delta^n\setminus\partial\Delta^n$ of $\Delta^n$; we use barycentric coordinates $u_0,\ldots, u_n$ on $V_n$, giving us the identification
\[
V_n=\Spec k[u_0,\ldots, u_n, (u_0\cdot\ldots\cdot u_n)^{-1}]/\sum_iu_i-1.
\]
We let $H\subset \P^n$ be the hyperplane $\sum_{i=1}^nX_i=X_0$ and let $1:=(1\:1\:\ldots\:1)\in\P^n(k)$. 

\begin{Def}\label{Def:PTCollapseMap} Let $F$ be finitely generated  field extension  of $k$ and let $w$ be a closed point of $V_{nF}$.  The {\em  Pontryagin-Thom collapse map} associated to $w$:
\[
PT_F(w):\Sigma^n_s\Spec F_+\to (\P^n_{F(w)}/H_{F(w)},1).
\]
 is the composition in $\sH_\bullet(k)$
\begin{multline*}
\Sigma^n_s\Spec F_+\xymatrix{\ar[r]^{\sigma_F}_\sim&} \Delta^n_F/\partial\Delta^n_F\xrightarrow{c_w} \Delta^n_F/(\Delta^n_F\setminus\{w\})
\xymatrix{\ar[r]^{mv_w} _\sim&}(\P^n_{F(w)}/H_{F(w)},1) 
\end{multline*}
 for specific choices of the isomorphisms in this composition, to be filled in below.
 \end{Def}  

The map $\sigma_F$ is the standard one given by the contractibility of $\Delta^n$ and all its faces, which gives an isomorphism in $\sH_\bullet(k)$ of $\Delta^n/\partial\Delta^n$ with the constant presheaf on the simplicial space  $\Delta_n/\partial\Delta_n$:
\[
\Delta_n([m]):=\Hom_\Delta([m],[n])
\]
and $\partial\Delta_n([m])\subset \Delta_n([m])$ the set of non-surjective maps $f:[m]\to [n]$. The isomorphism $\Sigma^nS^0\cong \Delta_n/\partial\Delta_n$ in $\sH_\bullet$ thus gives the isomorphism 
\[
\sigma:\Sigma^nS^0\to \Delta^n/\partial\Delta^n
\]
in $\sH_\bullet(k)$ and thereby gives rise to the isomorphism in $\sH_\bullet(k)$ 
\begin{equation}\label{eqn:sigma}
\sigma_F:\Sigma^n_s\Spec F_+=\Spec F_+\wedge\Sigma^nS^0\xrightarrow{\id\wedge\sigma} \Spec F_+\wedge\Delta^n/\partial\Delta^n= \Delta^n_F/\partial\Delta^n_F.
\end{equation}

The map $c_w$ is the quotient map. The isomorphism
\[
mv_w:\Delta^n_F/\Delta^n_F\setminus\{w\}\to (\P^n_{F(w)}/H_{F(w)},1)
\]
is the Morel-Voevodsky purity isomorphism. This map depends in general on the choice of an isomorphism
$\psi_w:m_w/m_w^2\to F(w)^n$, where $m_w\subset\sO_{\Delta^n_F,w}$ is the maximal ideal;
in addition, we need to make explicit the role of  the chosen  base-point 1. For this, we go through the construction of the purity isomorphism, giving the explicit choices which lead to a well-defined choice of isomorphism $mv_w$.

We give $V_n\times\A^1\times\Delta^n$ coordinates $u_0,\ldots, u_n, x, t_0,\ldots, t_n$, with the $u_i$ the barycentric coordinates on $V_n$, $x$ the standard coordinate on $\A^1$ and the $t_i$ the barycentric coordinates on $\Delta^n$.
Let 
\[
(X_0, X_1,\ldots, X_n):=(x, \frac{t_1-u_1}{u_0},\ldots, \frac{t_n-u_n}{u_0}).
\]

The construction of $mv_w$ uses the blow-up of $\A^1\times\Delta_F^n$ along $0\times w$
\[
\mu_w:\Bl_{0\times w}\A^1\times\Delta_F^n\to  \A^1\times\Delta_F^n.
\] 
Let $E_w\subset \Bl_{0\times w}\A^1\times\Delta_F^n$ be the exceptional divisor. Then $E_w$ is an $F(w)$-scheme. 

Suppose first that $w$ is separable over $F$. The closed point $0\times w$ of $\A^1\times \Delta^n_F$ has the canonical lifting to the closed point $0\times w$ of  $\A^1\times \Delta^n_{F(w)}$; let $m_{0\times w}\subset \sO_{\A^1\times \Delta^n_F,0\times w}$ and $m'_{0\times w}\subset \sO_{\A^1\times \Delta^n_{F(w)},0\times w}$ denote the respective maximal ideals. As $w$ is separable over $F$, the projection $p:\A^1\times \Delta^n_{F(w)}\to \A^1\times \Delta^n_F$ induces an isomorphism of graded $F(0\times w)$-algebras
\[
p^*:\oplus_{m\ge0} m_{0\times w}^m/m_{0\times w}^{m+1}\to \oplus_{m\ge0} m_{0\times w}^{\prime m}/m_{0\times w}^{\prime m+1}.
\]

The functions $(X_0, X_1(w),\ldots, X_n(w))$ give generators for the maximal ideal $m'_{0\times w}$; as
\[
E_w=\Proj_{F(0\times w)}\oplus_{m\ge0} m_{0\times w}^m/m_{0\times w}^{m+1}\cong \Proj_{F(0\times w)}\oplus_{m\ge0} m_{0\times w}^{\prime m}/m_{0\times w}^{\prime m+1}
\]
the image  $(x_0, x_1(w),\ldots, x_n(w))$ of $(X_0, X_1(w),\ldots, X_n(w))$ in $m'_{0\times w}/m_{0\times w}^2$ give  homogeneous coordinates for $E_w$, defining an isomorphism
\[
q_w:=(x_0\:x_1(w)\:\ldots\:x_n(w)):E_w\to \P^n_{F(w)}.
\]
Let  $H(w)\subset E_w$ be the pull-back of $H_{F(w)}$ via $q_w$, and let $1_w=q_w^{-1}(1)$. 

The proper transform $\mu_w^{-1}[\A^1\times w]\subset \Bl_{0\times w}\A^1\times\Delta_F^n$ maps isomorphically to $\A^1\times w$ via $\mu_w$, and intersects $E_w$ in a closed point $\bar{w}$ lying over $0\times w$.

\begin{lem}\label{lem:gen} 1. For all $w\in V_{nF}$, we have $1_w\neq \bar{w}$ and $\bar{w}\not\in H(w)$.\\
2. $q_w(\bar{w})=(1\:0\:\ldots\:0)$
\end{lem}

\begin{proof} Clearly (2) implies (1). For (2), $q_w(\bar{w})$ is the image of  $1\times w$ under 
\[
(X_0\: X_1(w)\:\ldots\:X_n(w)):\A^1\times\Delta^n_{F(w)}\setminus\{0\times w\}\to \P^n_{F(w)},
\]
which is $(1\:0\:\ldots\: 0)$.
\end{proof}

Additionally, the quotient map
\[
r_w:(\P^n_{F(w)}/H_{F(w)},1)\to \P^n_{F(w)}/(\P^n_{F(w)}\setminus \{(1\:0\:\ldots\:0)\})
\]
is an isomorphism in $\sH_\bullet(k)$, since projection from $(1\:0\:\ldots\:0)$ realizes $\P^n_{F(w)}\setminus\{(1\:0\:\ldots\:0)\}$ as an $\A^1$-bundle over $\P^{n-1}_{F(w)}$ with section $H_{F(w)}$.

This gives us the sequence of isomorphisms in $\sH_\bullet(k)$:
\[
E_w/(E_w\setminus \{\bar{w}\})\xrightarrow{q_w}\P^n_{F(w)}/(\P^n_{F(w)}\setminus\{(1\:0\:\ldots\:0)\})\xleftarrow{r_w}
(\P^n_{F(w)}/H_{F(w)},1).
\]

In case $w$ is not separable over $F$, we choose any set of parameters $X_1(w)$, $\ldots$, $X_n(w)$ for $m_w$  such that,  taking $X_0=x$, the isomorphism $E_w\to \P^n_w$ defined by the sequence $x_0, x_1(w),\ldots, x_n(w)$ satisfies the condition of lemma~\ref{lem:gen} ($F$ is infinite, so (1) is satisfied for a general choice; the condition (2) is satisfied for all choices). We then  proceed as above.

Morel-Voevodsky show that the inclusions $i_w:E_w\to \Bl_{0\times w}\A^1\times\Delta_F^n$ and
$i_1:\Delta_F^n=1\times\Delta^n_F\to \A^1\times\Delta^n_F$ induce  isomorphisms
\begin{align*}
&\bar{i}_w:E_w/(E_w\setminus \{\bar{w}\})\to \Bl_{0\times w}\A^1\times\Delta_F^n/(\Bl_{0\times w}\A^1\times\Delta_F^n\setminus \mu_w^{-1}[\A^1\times w])\\
&\bar{i}_1:\Delta^n_F/(\Delta^n_F\setminus\{w\})\to  \Bl_{0\times w}\A^1\times\Delta_F^n/(\Bl_{0\times w}\A^1\times\Delta_F^n\setminus \mu_w^{-1}[\A^1\times w])
\end{align*}
in $\sH_\bullet(k)$ (see the proof of  \cite[Theorem 3.2.23]{MorelVoev}).

\begin{Def}\label{Def:PurityIsomv_w} The purity isomorphism
\[
mv_w:\Delta^n_F/(\Delta^n_F\setminus\{w\})\xrightarrow{\sim} (\P^n_{F(w)}/H_{F(w)},1).
\] 
is defined as the composition
\begin{align*}
\Delta^n_F/(\Delta^n_F\setminus\{w\})&\xymatrix{\ar[r]^{\bar{i}_1}_\sim&}   \Bl_{0\times w}\A^1\times\Delta_{F}^n/(\Bl_{0\times w}\A^1\times\Delta_{F}^n\setminus 
 \mu_w^{-1}[\A^1\times w])\\
&\xymatrix{&\ar[l]_{\bar{i}_w}^\sim} E_w/(E_w\setminus \{\bar{w}\})\\
&\xymatrix{\ar[r]^{q_w}_\sim&} \P^n_{F(w)}/\P^n_{F(w)}\setminus\{(1\:0\:\ldots\:0)\}\\
&\xymatrix{&\ar[l]_{r_w}^\sim}(\P^n_{F(w)}/H_{1 F(w)},1).
\end{align*}
\end{Def}

In case $w$ is an $F$-rational point of $\Delta^n_F$, we have another description of $mv_w$.    The map
\[
q_w^{-1}\circ (X_0\:X_1(w)\:\ldots\:X_n(w)):\A^1\times\Delta^n_{F(w)}\setminus\{0\times w\}\to E_w
\]
extends to a morphism
\[
p_w:\Bl_{0\times w}\A^1\times\Delta_F^n\to E_w
\]
making $\Bl_{0\times w}\A^1\times\Delta_F^n$ an $\A^1$-bundle over $E_w$ with section $i_w$, and thus $p_w$ induces an isomorphism in $\sH_\bullet(k)$
\[
\bar{p}_w:\Bl_{0\times w}\A^1\times\Delta_F^n/(\Bl_{0\times w}\A^1\times\Delta_F^n\setminus 
 \mu_w^{-1}[\A^1\times w])\to 
E_w/(E_w\setminus \{\bar{w}\})
\]
inverse to  $\bar{i}_w$. Thus

\begin{lem}\label{lem:MVRatl} Suppose $w$ is in $\Delta^n(F)$. Then 
\[
mv_w=r_w^{-1}\circ q_w\circ\bar{p}_w\circ \bar{i}_1.
\]
\end{lem}

We can further simplify the above description of $mv_w$ by noting:
\begin{lem}  Suppose $w$ is in $\Delta^n(F)$. Let 
\begin{equation} \label{eqn:phi}
\phi_w:\Delta^n_F/(\Delta^n_F\setminus\{w\})\to  \P^n_F/\P^n_F\setminus\{(1\:0\:\ldots\:0)\}
\end{equation}
be the map induced by 
\[
(1\: X_1(w)\:\ldots\: X_n(w)):\Delta^n_F\to \P^n_F.
\]
Then  $\phi_w  =q_w\circ \bar{p}_w\circ   \circ \bar{i}_1$, hence
$mv_w=r_w^{-1}\circ\phi_w$.
\end{lem}

\begin{proof} The identity $mv_w=r_w^{-1}\circ\phi_w$
 follows directly from our description above of the maps $q_w$ and $\bar{p}_w$ and lemma~\ref{lem:MVRatl}.
 \end{proof}
 
 Altogether, this gives us the formula, for $w\in \Delta^n(F)$,
\begin{equation}\label{eqn:PT}
PT_F(w)= r_w^{-1}\circ \phi_w\circ c_w\circ \sigma_F.
\end{equation}

\section{$(\P^n/H,1)$ and $\Sigma^n_s\G_m^{\wedge n}$}\label{sec:Susp1}
Our main task in this section is to construct an explicit isomorphism
\[
\alpha:(\P^n/H,1) \xrightarrow{\sim} \Sigma^n_s\G_m^{\wedge n}.
\]

We first recall some elementary constructions involving homotopy colimits over subcategories of the $n$-cube. Let $\sC$ be a small category and let $\sF:\sC\to \Spc(k)$ be a functor. Let $\sN:\Delta^\op\to\Sets$ be the nerve of $\sC$. For 
\[
\sigma=(s_0\xrightarrow{f_1}s_1\to\ldots\xrightarrow{f_n}s_n)\in\sN_n(\sC)
\]
define $\sF(\sigma):=\sF(s_0)$.  Bousfield-Kan \cite{BousfieldKan}  define $\underline{\hocolim}\,\sF$ to be  the simplicial object of $\Spc(k)$ with $n$-simplices
\[
\underline{\hocolim}\,\sF_n:=\amalg_{\sigma\in\sN_n(\sC)}\sF(\sigma);
\]
for $g:[n]\to [m]$ in $\Delta$, 
\[
\underline{\hocolim}\sF(g):\underline{\hocolim}\,\sF_m\to \underline{\hocolim}\,\sF_n
\]
is the map sending $(\sF(s_0), \sigma=(s_0,\ldots, s_m))$ to $(\sF(s'_0), \sigma'=(s'_0,\ldots, s'_n))$, with $\sigma'=\sN(g)(\sigma)$, $s'_0=s_{g(0)}$ and the map $\sF(s_0)\to \sF(s'_0)$ is $\sF(s_0\to s_{g(0)})$.  $\hocolim\sF$ is the geometric realization of  $\underline{\hocolim}\,\sF$.

For a functor $\sF:\sC\to \Spc_\bullet(k)$ we use essentially the same definition of $\underline{\hocolim}\,\sF$ as a simplicial object of $ \Spc_\bullet(k)$, replacing disjoint union $\amalg$ with pointed union $\vee$, and we use the pointed version of geometric realization to define 
$\hocolim\sF$ in $\Spc_\bullet(k)$. Concretely, $\hocolim\sF$ is the co-equalizer of 
\[
\vee_{g:[n]\to [m]} \underline{\hocolim}\,\sF_m\wedge\Delta^n_+
\xymatrix{\ar@<3pt>[r]\ar@<-3pt>[r]&}\vee_n\underline{\hocolim}\,\sF_n\wedge\Delta^n_+
\]

The essential property of $\hocolim$ we will need is the following:
\begin{prop}[\hbox{\cite{BousfieldKan}}] Let $\sC$ be a finite category, $\sF, \sG:\sC\to \Spc_\bullet(k)$ functors, and $\vartheta:\sF\to\sG$ a natural transformation. Suppose that $\vartheta(c):\sF(c)\to \sG(c)$ is an isomorphism in $\sH_\bullet(k)$ for each $c\in\sC$. Then
\[
\hocolim\vartheta:\hocolim\sF\to \hocolim\sG
\]
is an isomorphism in $\sH_\bullet(k)$. The analogous result holds after replacing $\Spc_\bullet(k)$ and $\sH_\bullet(k)$ with $\Spc(k)$ and $\sH(k)$.
\end{prop}
This is of course just a special case of the general result valid for functors from a small (not just finite) category to a proper simplicial model category. See for example \cite{Jardine} for details.

\begin{rem}\label{rem:Bisimp} Let $\sF:\sC\to \Spc(k)$ be a functor. Suppose our index category $\sC$ is a product $\sC_1\times\sC_2$. We may form the bi-simplicial object  $\underline{\hocolim}^2\,\sF$ of $\Spc(k)$, with $(n,m)$-simplices 
\[
\underline{\hocolim}^2\,\sF_{n,m}:=\amalg_{(\sigma_1,\sigma_2)\in\sN(C_1)_n\times\sN(C_2)_m}\sF(\sigma_1,\sigma_2)
\]
where $\sF(\sigma_1,\sigma_2)=\sF(s_0\times s'_0)$ if $\sigma=(s_0\to\ldots)$ and $\sigma'=(s'_0\to\ldots)$; the morphisms are defined similarly.

As $\sN(\sC_1\times\sC_2)$ is the diagonal simplicial set associated to the bi-simplicial set $\sN(\sC_1)\times\sN(\sC_2)$, it follows that $\underline{\hocolim}\,\sF$ is the diagonal simplicial object of $\Spc(k)$ associated to the bi-simplicial object $\underline{\hocolim}^2\,\sF$, and thus we have the natural isomorphism of geometric realizations
\[
\hocolim\sF=|\underline{\hocolim}\,\sF|\cong |\underline{\hocolim}^2\,\sF|
\]
in $\Spc(k)$. Similar remarks hold in the pointed case.
\end{rem}

Let $\square^{n+1}$ be the poset of subsets of  $[n]$, ordered under inclusion. For $J\subset J'\subset\{0,\ldots,n\}$, we let $\square^{n+1}_{J\le *\le J'}$ be the full subcategory of subsets $I$ with $J\subset I\subset J'$, $\square^{n+1}_{J< *\le J'}$ the full subcategory of subsets $I$ with $J\subsetneq I\subset J'$, etc. We sometimes omit $J$ if $J=\0$ or $J'$ if $J'=[n]$.

If $|J|=r+1$, we let $i'_J:\{0,\ldots, r\}\to J$ be the unique order-preserving bijection, and let
$i_J:\square^{r+1}\to \square^{n+1}_{*\le J}$
be the resulting isomorphism of categories. Clearly $i_J$ induces the isomorphism of subcategories
$i_J:\square^{r+1}_{*<[r]}\to \square^{n+1}_{*<J}$.

We will be using the following elementary constructions. Let $\sF:\square^{r+1}_{<[r]}\to \Spc_\bullet(k)$ be a functor and take $n>r$. Identifying $\square^{r+1}_{*<[r]}$ with $\square^{n+1}_{*<[r]}$ via the inclusion $[r]\subset [n]$, extend $\sF$ to a functor
\[
\sigma^{n-r}\sF:\square^{n+1}_{*<[n]}\to \Spc_\bullet(k)
\]
by setting $\sigma^{n-r}\sF(J)=*$ if $J$ is not a proper subset of $[r]$. 

Similarly, let $\sG:\square^{r+1}_{*<[r]}\times\square^s\to \Spc_\bullet(k)$ be a functor and take $n>r$. Identifying
$\square^{r+1}_{*<[r]}\times\square^s$ with a full subcategory of $\square^{r+1}_{*<[r]}\times\square^{n-r+s}$ via the inclusion $[s]\subset [n-r+s]$, extend $\sG$ to a functor
\[
c^{n-r}\sG:\square^{r+1}_{*<[r]}\times\square^{n-r+s}\to \Spc_\bullet(k)
\]
by setting $c^{n-r}\sG(J,I)=*$  if $I\not\subset[s]$.

\begin{ex} Let $\sX$ be in $\Spc_\bullet(k)$. Noting that $\square^1_{*<[0]}$ is the one-point category, we  write $\sX$ for the functor $\square^1_{*<[0]}\to \Spc_\bullet(k)$ with value $\sX$. This gives us the functors 
\begin{align*}
&c^n\sX:\square^{n+1}\to \Spc_\bullet(k),\\
&\sigma^n\sX:\square^{n+1}_{*<[n]}\to \Spc_\bullet(k).
\end{align*}
Explicitly, $c^n\sX(\0)=\sigma^n\sX(\0)=\sX$ and both functors have value $*$ at $J\neq\0$.
\end{ex}

\begin{lem} \label{lem:Susp}   There are natural isomorphisms
\begin{align*}
&\Pi_c:\hocolim c^{n-r}\sG \to \hocolim \sG\wedge([0,1],1)^{\wedge n-r}\\
&\Pi_\sigma:\hocolim\sigma^{n-r}\sF \to \Sigma_s^{n-r}\hocolim\sF.
\end{align*}
in $\Spc_\bullet(k)$
\end{lem}

\begin{proof} We proceed by induction on $n-r$; it suffices to handle the case $r=n-1$. We first take care of the isomorphism $\Pi_c$.

Via remark~\ref{rem:Bisimp}, it suffices to give an isomorphism 
\[
|\underline{\hocolim}^2\, c^1\sG|\cong \hocolim\sG\wedge([0,1],1),
\]
where we use the product decomposition $\square^{r+1}_{*<[r]}\times\square^{s+1}=(\square^{r+1}_{*<[r]}\times\square^s)\times\square^1$.
Fix an $m$ simplex $\sigma$ of $\sN(\square^{r+1}_{*<[r]}\times\square^s)$ and let 
\[
c\sG_\sigma:\square^1\to \Spc_\bullet(k)
\]
be the functor $c\sG_\sigma(\0)=\sG(\sigma)\wedge\Delta^m_+$, $c\sG_\sigma([0])=*$. Then
\[
\hocolim c\sG_\sigma\cong \sG(\sigma)\wedge\Delta^m_+\wedge ([0,1],1)
\]
with the isomorphism natural in $\sigma$. The result follows directly from this.

Next, given  $\sF:\square^{n}_{<[n-1]}\to \Spc_\bullet(k)$, let $c'\sF:\square^n\to \Spc_\bullet(k)$ be the extension of $\sF$ to $\square^n$ defined by setting $c'\sF([n-1])=*$. We claim there is a natural isomorphism
\[
\hocolim c'\sF\cong \hocolim\sF\wedge ([0,1],1)
\]
in $\Spc_\bullet(k)$.

Indeed, we have the bijection (for $m>0$)
\[
\sN(\square^n)_m=\sN(\square^n_{*<[n-1]})_m\amalg\sN(\square^n_{*<[n-1]})_{m-1}
\]
with the first component coming from the inclusion of $\square^n_{*<[n-1]}$ in $\square^n$, and the second arising by sending $\sigma=(s_0\to\ldots\to s_{m-1})$ to $(s_0\to\ldots\to s_{m-1}\to[n-1])$. For $m=0$, the same construction gives 
\[
\sN(\square^n)_0=\sN(\square^n_{*<[n-1]})_0\amalg\{[n-1]\}.
\]
As, for a simplicial set $C$, the $m$-simplices of $C\times[0,1]/C\times 1$ have exactly the same description, our claim follows easily.

Finally, we can write the category $\square^{n+1}_{*<[n]}$ as a (strict) pushout
\[
\square^{n+1}_{*<[n]}=\square^{n}\amalg_{\square^n_{*<[n]}} \square^n_{*<[n]}\times \square^1.
\]
This leads to an isomorphism of $\hocolim\sigma^1\sF$ as a pushout
\begin{align*}
\hocolim\sigma^1\sF&\cong \hocolim c'\sF\vee_{\hocolim \sF} \hocolim c\sF\\
&\cong \hocolim \sF\wedge  ([0,1],1)\vee_{\hocolim \sF\wedge0_+}  \hocolim \sF\wedge  ([0,1],1)\\
&=\Sigma^1_s\hocolim \sF.
\end{align*}
\end{proof}

As in section~\ref{sec:PTCollapse}, let $H\subset \P^n$ be the hyperplane $\sum_{i=1}^nX_i=X_0$ and let $1:=(1\:1\:\ldots\: 1)\in\P^n(k)$. We define an isomorphism $\alpha: (\P^n/H,1)\xrightarrow{\sim} \Sigma_s^n\G_m^{\wedge n}$
in $\sH_\bullet(k)$ as follows:

Let $U_i\subset \P^n$ be the standard affine open subset $X_i\neq 0$. We identify $U_i$ with $\A^n$ in the usual way via coordinates $(X_0/X_i,\ldots\hat{X_i/X_i},\ldots, X_n/X_i)$, which we write as $x^i_1,\ldots, x^i_n$, or simply $x_1,\ldots, x_n$. For each index set $I\subset \{0,\ldots, n\}$, we have the intersection
\[
U_I:=\cap_{i\in I}U_i.
\]
For $I=\{i_1<\ldots<i_r\}$, we use coordinates in $U_{i_1}$ to identify 
\[
U_I\cong\Spec k[x_1,\ldots, x_n, x_{i_2}^{-1},\ldots, x_{i_r}^{-1}]\cong \A^{n-|I|+1}\times\G_m^{|I|-1}.
\] 

The open cover $\sU:=\{U_0,\ldots, U_n\}$ of $\P^n$ identifies $\P^n$ (in $\sH(k)$) with the homotopy colimit over $\square^{n+1}_{*<[n]}$ of the functor
\begin{align*}
\sP^n_{\sU}:\square^{n+1}_{*<[n]}&\to \Spc(k)\\
\sP^n_{\sU}(J)&:= U_{J^c}.
\end{align*}
We thus have the functor
\begin{align*}
\sP^n_{\sU,1}:\square^{n+1}_{*<[n]}&\to \Spc_\bullet(k)\\
\sP^n_{\sU,1}(J)&:= (U_{J^c},1)
\end{align*}
and the isomorphism in $\sH_\bullet(k)$, $\hocolim\sP^n_{\sU,1}\cong (\P^n,1)$.

 Next, we note that the hyperplane $H\subset \P^n$ is covered by the affine open subsets $U_1,\ldots, U_n$. The open cover $\sU_1:=\{H\cap U_1,\ldots, H\cap U_n\}$ of $H$  identifies $H$ (in $\sH(k)$) with the homotopy colimit over $\square^{n+1}_{*<[n]}$ of the functor
\begin{align*}
\sH_{\sU_1}:&\square^{n+1}_{*<[n]}\to \Spc(k)\\
\sH_{\sU_1}(J)&:= \begin{cases} H\cap U_{J^c}&\text{ for }0\in J\\ \0&\text{ for }0\not\in J.\end{cases}
\end{align*}

Let 
 \[
\sP^n_{\sU,1}/\sH_{\sU_1}:\square^{n+1}_{*<[n]}\to \Spc_\bullet(k)
\]
be the functor defined by
 \[
\sP^n_{\sU,1}/\sH_{\sU_1}(J):=\begin{cases} (U_{J^c}/H\cap U_{J^c},1)&\text{ for } 0\in J\\ (U_{J^c},1)&\text{ for }0\not\in J.\end{cases}
\]
By our discussion, the maps $\sP^n_{\sU,1}/\sH_{\sU_1}(J)\to (\P^n/H, 1)$ induced by the inclusions $U_{J^c}\hookrightarrow \P^n$
give rise to an isomorphism in $\sH_\bullet(k)$
\[
\epsilon_1:\hocolim  \sP^n_{\sU,1}/\sH_{\sU_1}\to (\P^n/H, 1).
\]
To simplify the notation, we denote $\sP^n_{\sU,1}/\sH_{\sU_1}$ by $\sF$ for the next few paragraphs.

We claim that, for each $J\neq\0$ with $0\not\in J$, we have 
$(U_{J}/(H\cap U_{J}),1)\cong *$
in $\sH_\bullet(k)$. Indeed, suppose for example that $n\in J$, and use coordinates $(x^n_1,\ldots, x^n_n)=(X_0/X_n,\ldots, X_{n-1}/X_n)$ on $U_J$.  We have the projection
\begin{align*}
&p:U_J\to U_J\cap(X_0=0)\\
&p(x^n_1,\ldots, x^n_n)=(0,x^n_2,\ldots, x^n_n).
\end{align*}
Since $0\not\in J$, $x^n_1$ is not inverted on $U_J$, and thus $p$ makes $U_J$ an $\A^1$-bundle over 
$U_J\cap(X_0=0)$. $p$ has the section 
\[
s(0,x^n_2,\ldots, x^n_n):=(1+\sum_{i=2}^nx^n_i,x^n_2,\ldots, x^n_n),
\]
identifying $U_J\cap(X_0=0)$ with $H\cap U_J$; this together with homotopy invariance in $\sH_\bullet(k)$ proves our claim. Thus $\sF(J)\cong *$ in $\sH_\bullet(k)$ for all $J$ with $0\in J$. In addition $\sF(\{1,2,\ldots, n\})=(U_0,*)\cong (\A^n,*)$, which is also isomorphic to $*$ in $\sH_\bullet(k)$.

Let $i_0:\square^n_{*<[n-1]}\to \square^{n+1}_{*<[n]}$  inclusion functor induced by the inclusion $[n-1]\to [n]$ sending $i\in[n-1]$ to $i+1$, and let $\omega:\square^{n+1}_{*<[n]}\to \square^{n+1}_{*<[n]}$ be the automorphism induced by the cyclic permutation $\omega$ of $[n]$,
\[
\omega(i):=\begin{cases}i+1&\text{ for }0\le i<n\\0&\text{ for }i=n.\end{cases} 
\]

Let
\[
\sF_{|0}:\square^n_{*<[n-1]}\to \Spc_\bullet(k)
\]
be the functor  $\sF\circ i_0$. We have the evident quotient map
$q:\sF\circ\omega\to \sigma^1\sF_{|0}$, 
which by our discussion above is a term-wise isomorphism in $\sH_\bullet(k)$. By lemma~\ref{lem:Susp}, $q$ induces the isomorphisms in $\sH_\bullet(k)$
\begin{equation}\label{eqn:SuspIso}
\hocolim\sF\to \hocolim\sigma^1\sF_{|0}\to \Sigma^1_s\hocolim\sF_{|0}.
\end{equation}

We now turn to the functor $\sF_{|0}$. This is just the punctured $n$-cube corresponding to the open cover $\sU':=\{U_0\cap U_1,\ldots, U_0\cap U_n\}$ of $U_0\setminus (1\:0\:\ldots\:0)$ (with base-point 1), i.e. $(\A^n\setminus 0,1)$. We thus have the isomorphism in $\sH_\bullet(k)$
\[
\hocolim\sF_{|0}\cong (U_0\setminus (1\:0\:\ldots\: 0),1)\cong (\A^n\setminus 0,1).
\]
Let $C\subset U_0\setminus (1\:0\:\ldots\: 0)$ be the union of the affine hyperplanes $x_i^0=1$, $i=1,\ldots, n$. As the inclusion $1\to C$ is an isomorphism in $\sH(k)$, we have the isomorphism in $\sH_\bullet(k)$
\[
 (U_0\setminus (1\:0\:\ldots\: 0),1)\cong  U_0\setminus (1\:0\:\ldots\: 0)/C.
\]
Letting $\bar{\sF}_{|0}$ be the quotient of $\sF_{|0}$ given by
\[
\bar{\sF}_{|0}(J)=U_{i_0(J)^c}/C\cap U_{i_0(J)^c},
\]
we thus have the isomorphisms in $\sH_\bullet(k)$
\[
\hocolim\sF_{|0}\cong \hocolim\bar{\sF}_{|0}\cong U_0\setminus (1\:0\:\ldots\: 0)/C.
 \]
 
On the other hand, for each $J\subsetneq\{1,\ldots, n\}$, the inclusion $C\cap U_0\cap U_J\to U_0\cap U_J$ is an isomorphism in 
 $\sH(k)$, and thus $\bar{\sF}_{|0}(J)\cong*$ for all $J\neq\0$. Since $\bar{\sF}_{|0}(\0)\cong \G_m^{\wedge n}$ we have the quotient map
 $ \bar{\sF}_{|0}\to \sigma^{n-1}\G_m^{\wedge n}$;
 our discussion together with lemma~\ref{lem:Susp} thus gives us the isomorphism in $\sH_\bullet(k)$
 \[
 \hocolim\sF_{|0}\cong \hocolim\bar{\sF}_{|0}\cong \Sigma^{n-1}_s\G_m^{\wedge n}.
 \]
Together with \eqref{eqn:SuspIso}, this gives us the sequence of isomorphisms in $\sH_\bullet(k)$
\[
 (\P^n/H,1)\cong \hocolim\sP^n_{\sU,1}/\sH_{\sU_1}\cong \Sigma_s\hocolim\sF_{|0}\cong \Sigma^{n}_s\G_m^{\wedge n}.
\]
We denote the composition by
\begin{equation}\label{eqn:SuspIso1}
\alpha: (\P^n/H,1)\xrightarrow{\sim}\Sigma^n_s\G_m^{\wedge n}.
\end{equation}

Now that we have defined $\alpha$, we can complete our definition of the purity isomorphism  \eqref{eqn:MV}:
\begin{equation}\label{eqn:MVDef}
MV_w:=(\id_{w_+}\wedge\alpha)\circ mv_v
\end{equation}
(see definition~\ref{Def:PurityIsomv_w} for the definition of  $mv_w$).

\begin{rem} Take $n>1$. Let $H_\infty\subset\P^n$ be the hyperplane $X_0=0$ and for  let $C_1\subset U_1$ be the union of the hyperplanes $x^1_i=1$, $i=1,\ldots, n$.   Let  $\sG^\infty_n:\square^{n+1}_{*<[n]}\to\Spc_\bullet(k)$ be the functor
\[
\sG^\infty_n(J):=\begin{cases} U_{J^c}/H_\infty \cap U_{J^c}&\text{ for } 1\in J\\
U_{J^c}/[(H_\infty\cup C_1) \cap U_{J^c}]&\text{ for } 1\not\in J.
\end{cases}
\]
We note that  the inclusion $(0\:1\:0\:\ldots\: 0)\to H_\infty\cap C_1$ is an $\A^1$-weak equivalence; using this it is easy to modify the arguments used in this section to show that the identity map $\sG^\infty(\0)\to \G_m^{\wedge n}$ extends to a map of functors
$\sG^\infty_n\to \sigma^n\G_m^{\wedge n}$, 
which is a termwise isomorphism in $\sH_\bullet(k)$, giving us the isomorphism
\[
\hocolim\sG^\infty_n\cong \Sigma^n_s\G_m^{\wedge n}
\]
 in $\sH_\bullet(k)$. Furthermore, we have the sequence of isomorphisms in $\sH_\bullet(k)$:
 \[
 \P^n/H_\infty\to \P^n/[H_\infty\amalg_{C_1\cap H_\infty\cap U_1} C_1]\to \hocolim\sG^\infty_n.
 \]
 Putting these together gives us the isomorphism
\begin{equation}\label{eqn:Iso2}
 \alpha_\infty: \P^n/H_\infty \to \Sigma^n_s\G_m^{\wedge n}
\end{equation}
 in $\sH_\bullet(k)$. 
 
 For $n=1$, we note that $H=1$, so $(\P^1/H,1)=(\P^1, H)$. To define $\alpha_\infty$, we just compose $\alpha:(\P^1/H,1)\to \Sigma_s\G_m$ with the isomorphism $\tau:(\P^1, H_\infty)\to (\P^1,H)$ given by
 \[
 \tau(X_0\:X_1)=(X_1-X_0\:X_1).
 \]
We will use these models for $\Sigma^n_s\G_m^{\wedge n}$ to construct transfer maps in \S\ref{sec:Transfer}.
\end{rem}

\section{The suspension of a symbol}\label{sec:Susp2}
Let $\tilde\rho:V_n\to\G_m^n$ be the map
\[
\rho(u_0,\ldots, u_n):=(-\frac{u_1}{u_0},\ldots, -\frac{u_n}{u_0}).
\]
Composing with the quotient map $\G_m^n\to\G_m^{\wedge n}$ gives us the map
$\rho:V_{n+}\to \G_m^{\wedge n}$.
Our next main task is to give  an explicit algebro-geometric description of $\Sigma^n_s\rho$. More generally, for
$f:T\to V_n$
a morphism in $\Sm/k$, we will give a description of $\Sigma^n_s(\rho\circ f)$. We begin by giving a description of $\Sigma^n_sT_+$ as a certain homotopy colimit.

For this,   consider the scheme $\A^1\times\Delta^n$, with coordinates $x, t_0,\ldots, t_n$:
\[
\A^1\times\Delta^n=\Spec k[x,t_0,\ldots, t_n]/\sum_it_i-1.
\]
For $i=1,\ldots, n$, let $U'_i\subset \A^1\times\Delta^n$ be the subscheme defined by $t_i=0$, and let 
$U'_0\subset \A^1\times\Delta^n$ be the subscheme defined by $x=1$.  For $I\subset\{0,\ldots, n\}$, let
$U'_I:=\cap_{i\in I}U'_i$,
the intersection taking place in $ \A^1\times\Delta^n$. This gives us the punctured $n+1$-cube
\[
\hat{\sG}_n^T:\square^{n+1}_{*<[n]}\to \Spc(k)
\]
with $\hat{\sG}^T_n(J):=T\times  U'_{J^c}$.

As above, use barycentric coordinates $u_0,\ldots, u_n$ for $V_n$. We pull these back to $T$ via $f$, and write $u_i$ for $f^*(u_i)$, letting the context make the meaning clear. Set
\[
(X_0, X_1,\ldots, X_n):=(x, \frac{t_1-u_1}{u_0},\ldots, \frac{t_n-u_n}{u_0})
\]
and set
\[
(x^i_1,\ldots, x^i_n):=(X_0/X_i,\ldots,\widehat{X_i/X_i},\ldots, X_n/X_i);\quad i=0,\ldots, n.
\]
Inside $T\times\A^1\times\Delta^n$, we have the ``hyperplane" $H(T)$ defined by
\[
\sum_{i=1}^nX_i=X_0.
\]

Fix an index $I=(i_0,\ldots, i_r)$ with $0\le i_0<\ldots<i_r\le n$, and write the complement of $I$ in $\{0,\ldots, n\}$ as $I^c=(j_1,\ldots, j_{n-r})$ with $j_1<\ldots<j_{n-r}$. We have the isomorphism
\[
\phi_I:=\id\times(x^{i_0}_{j_1},\ldots, x^{i_0}_{j_{n-r}}):T\times U'_I\to T\times\A^{n-r}.
\]
In addition, let $H_I\subset \A^{n-r}$ be the hyperplane defined by 
\[\sum_{\ell =1}^{n-r}x_\ell=1\text{ if } i_0=0,\quad
\sum_{\ell =2}^{n-r}x_\ell=x_1\text{ if } i_0>0.
\]
Then  $\phi_I$ restricts to an isomorphism of $H(T)\cap T\times U'_I$ with $T\times H_I$, and thus the projection
$p_1:H(T)\cap T\times U'_I\to T$
and inclusion $\iota:H(T)\cap T\times U'_{I}\to T\times U'_{I}$
are isomorphisms in $\sH(k)$.

For $J\subsetneq[n]$, $J\neq \0$, define
$\sG^{T\prime}_n(J)$ to be the pushout in the diagram
\[
\xymatrix{
H(T)\cap T\times U'_{J^c}\ar@{^(->}[r]^-\iota\ar[d]_{p_1}&\hat{\sG}^T_n(J)\ar@{=}[r]\ar@{-->}[d]^{i(J)}&T\times U'_{J^c}\\
T\ar@{-->}[r]_-{s_J}&\sG^{T\prime}_n(J).
}
\]
Since $\iota$ is a cofibration and a weak equivalence in $\Spc(k)$, so is  $s_J$. As  $p_1$ is also a weak equivalence  in $\Spc(k)$, $i(J)$ 
is a weak equivalence in $\Spc(k)$ as well.

We set 
\[
\sG_n^{T\prime}(\0):=\hat\sG_n^{T}(\0)=T\times U_{[n]}'\cong T.
\]
This defines for us the functor
\[
\sG^{T\prime}_n:\square^{n+1}_{*<[n]}\to \Spc(k)
\]
that fits into a diagram ($T$ the constant functor)
\[
\xymatrix{
&\hat{\sG}^{T}_n\ar[d]^i\\
T\ar[r]_-s&\sG^{T\prime}_n
}
\]
with $i$ and $s$ term-wise isomorphisms in $\sH(k)$ and $s$ a term-wise cofibration in $\Spc(k)$.

For $n=1$,  define
\[
\sG_1^T(J):=\begin{cases}\sG^{T\prime}_1(J)/s(T)&\text{ for }J\neq\0\\ \sG^{T\prime}_1(\0)_+\cong T_+&\text{ for }J=\0.\end{cases}
\]
giving us the functor
\[
\sG^T_1:\square^{2}_{*<[1]}\to \Spc_\bullet(k)
\]

For $n>1$, take $\0\neq J\subset[n]$ and let $\Pi'_J\subset\P^n$ be the dimension $n-|J|$ linear subspace defined by $\cap_{j\in J}(X_j=0)$. Let $\Pi_J\subset\P^n$ be the dimension $n-|J|+1$ linear space spanned by $1$ and $\Pi'_J$ and let $\A_J\subset \Pi_J$ be the affine space $\Pi_J\setminus \Pi'_J$. Since $\Pi'_J$ is not contained in $H$, the intersection $\A_J\cap H$ is a codimension one affine space $\A_{J,H}$ in $\A_J$. Clearly $\A_J\supset\A_{J'}$ for $J\subset J'$, so we have the functor
\begin{align*}
\A/\A_H:\square^{n+1}_{*<[n]}&\to \Spc(k)\\
J&\mapsto \A_{J^c}/\A_{J^c,H}.
\end{align*}

Let $*_J$ be the base-point in $\A_J/\A_{J,H}$ and let $s'_J:T\to T\times \A_J/\A_{J,H}$
be the morphism identifying $T$ with $T\times *_J$.  Let $1_J$ be the image of $1\in \A_J$ in the quotient $\A_J/\A_{J,H}$. We have the morphism $s'_{J,1}:T\to T\times\A_J/\A_{J,H}$
identifying $T$ with $T\times 1_J$. For $J\neq\0$, let  $\sG^T_n(J)$ be the push-out in the diagram
\[
\xymatrix{
T\amalg T\ar[r]^-{s'_{J^c}\times s'_{J^c,1}}\ar[d]_{s_J\amalg p}&T\times \A_{J^c}/\A_{J^c,H}\ar@{-->}[d]\\
\sG^{T\prime}_n(J)\amalg *\ar@{-->}[r]&\sG^T_n(J).
}
\]
where $p:T\to *$ is the canonical map; we give $\sG^T_n(J)$ the base-point $*$.  We set $\sG^T_n(\0)=T_+$ with its canonical base-point.  Using the functoriality of $\sG_n^{T\prime}$ and $\A/\A_H$ defines the functor
\begin{equation}\label{eqn:GTDiag}
\sG_n^T:\square^{n+1}_{*<[n]}\to \Spc_\bullet(k).
\end{equation}

\begin{lem}\label{lem:Contract} For each $J\neq\0$, $\sG_n^T(J)\cong*$ in $\sH_\bullet(k)$.
\end{lem}

\begin{proof}Take $J\subset[n]$, $J\neq\0$.  For $n=1$,   $s:T\to \sG^{T\prime}_1(J)$ is a cofibration and weak equivalence in $\Spc(k)$, and thus the quotient $ \sG^{T\prime}_1(J)/T$ is contractible.

For $n>1$, the morphisms  $s_J:T\to \sG^{T\prime}_n(J)$, $s'_J:T\to T\times \A_J/\A_{J,H}$ and $s'_{J,1}:T\to T\times \A_J/\A_{J,H}$ are cofibrations and weak equivalences in $\Spc(k)$; since $1_J\not\in \A_{J,H}$, the map
\[
s'_J\times s'_{J,1}:T\amalg T\to T\times \A_J/\A_{J,H}
\]
is a cofibration.

Let $\sG^{T\prime\prime}_n(J)$ be the push-out in the diagram
\[
\xymatrix{
T\ar[r]^-{s'_{J^c}}\ar[d]_{s_J}&T\times \A_{J^c}/\A_{J^c,H}\ar@{-->}[d]^\iota\\
\sG^{T\prime}_n(J)\ar@{-->}[r]&\sG^{T\prime\prime}_n(J).
}
\]
Then  $\iota$ is a cofibration and a weak equivalence, hence the same is true for the composition
\[
T\xrightarrow{s'_{J,1}}T\times \A_J/\A_{J,H}\xrightarrow{\iota} \sG^{T\prime\prime}_n(J).
\]
As $\sG_n^T(J)=\sG^{T\prime\prime}_n(J)/T$, 
it follows that $\sG_n^T(J)$ is contractible.
\end{proof}

Letting $\sT:\square^1_{*<[0]}\to \Spc_\bullet(k)$ be the functor $\sT(\0)=T_+$, we have the evident quotient map
$\sG^T_{n}\to\sigma^n\sT$, i.e., we send $\sG_{n}(\0)=T_+$ to $\sigma^n\sT(\0)=T_+$ by the identity map, and the other maps are the canonical ones $\sG^T_{n}(I)\to *$. 

By lemma~\ref{lem:Susp} and  lemma~\ref{lem:Contract}, this map induces an isomorphism
\begin{equation}\label{eqn:beta}
\beta^T:\hocolim\sG^T_n\to \Sigma^n_sT_+
\end{equation}
in $\sH_\bullet(k)$.

\begin{rem} The  functors $\sG^T_n$, $\sG^{T\prime}_n$ and $\hat{\sG}^T_n$ are all functors in $T$, where for example $g:T'\to T$ gives the morphism $\hat{\sG}_n(f):\hat{\sG}_n^{T'}\to \hat{\sG}_n^T$ by the collection of maps
\[
f\times\id:T'\times U'_{J^c}\to T\times U'_{J^c}.
\]
The   map $\sG^{T\prime}_n\to \sG_n^T$ is natural in $T$, as is the map $\beta^T$.
\end{rem}

Let  $\Delta(V_n)\subset V_n\times\Delta^n$
be the graph of the inclusion $V_n\to\Delta^n$; by a slight abuse of notation, we write
$0\times \Delta(V_n)\subset V_n\times\A^1\times\Delta^n$
for the image of $0\times \Delta(V_n)\subset \A^1\times V_n\times\Delta^n$ under the exchange of factors
$\A^1\times V_n\times\Delta^n\to V_n\times\A^1\times\Delta^n$.

Define the morphism $\phi:V_n\times\A^1\times\Delta^n\setminus 0\times \Delta(V_n)\to \P^n$
by
\[
\phi(u_0,\ldots, u_n,x, t_0,\ldots, t_n):=(X_0\:X_1\:\ldots\:X_n),
\]
where as above $X_0=x$, $X_i=(t_i-u_i)/u_0$, $i=1,\ldots, n$.

Since  $V_n\times U'_i\cap 0\times \Delta(V_n)=\0$
for each $i=0,\ldots, n$, the restriction of $\phi$ to $\cup_{i=0}^nV_n\times U'_i$ is thus a morphism, and therefore gives a well-defined morphism of functors $\square^{n+1}_{*<[n]}\to \Spc(k)$, $\tilde\phi_*:\hat{\sG}_n^{V_n}\to \P^n$,
where $\P^n$ is the constant functor. 

Given a morphism $f:T\to V_n$, we compose $\tilde\phi_J$ with $f\times\id$, giving the morphism of functors
$\tilde\phi^T_*:\hat{\sG}_n^{T}\to \P^n$. 
Adjoining the projections $T\times U'_{J^c}\to T$ gives us the morphism of functors
$(p_1,\tilde\phi^T_*):\hat{\sG}_n^{T}\to T\times \P^n$.
Passing to the quotients, $(p_1,\tilde\phi^T_*)$ induces the map of functors
$(p_1,\phi^{T\prime}_*):\sG^{\prime T}_n\to T\times(\P^n/H)$.

We extend $(p_1,\phi^{T\prime}_*)$ to a map of functors $\square^{n+1}_{*<[n]}\to \Spc_\bullet(k)$
\[
p_1\wedge\phi^T_*:\sG^{T}_n\to T_+\wedge (\P^n/H,1)
\]
by using the inclusions $\A_{J^c}\to \P^n$, and sending the base-point in $T_+$ to the base-point in $T_+\wedge (\P^n/H,1)$. This gives us the map in $\Spc_\bullet(k)$
\begin{equation}\label{eqn:**}
\Phi^T:\hocolim\sG^T_n\to T_+\wedge (\P^n/H,1).
\end{equation}

\begin{lem}\label{lem:HauptLem} Let $f:T\to V_n$ be a morphism in $\Sm/k$. Then the diagram
\[
\xymatrixcolsep{30pt}
\xymatrix{
\Sigma_s^nT_+\ar[r]^-{\Sigma_s^n(\id_{T+}\wedge \rho\circ f)}&T_+\wedge\Sigma^n_s\G_m^{\wedge n}\\
\hocolim\sG^T_n\ar[r]_-{\Phi^T}\ar[u]^{\beta^T}&T_+\wedge (\P^n/H,1)\ar[u]_{\id\wedge\alpha}
}
\]
commutes in $\sH_\bullet(k)$.
\end{lem}

\begin{proof}  We work through our description of $\alpha$ and $\beta^T$, adding some intermediate steps. 

We introduce an additional functor
\begin{align*}
(\sP^n/\sH_{\sU},1):&\square^{n+1}_{*<[n]}\to \Spc_\bullet(k)\\
&J\mapsto (U_{J^c}/H\cap U_{J^c},1)
\end{align*}
By Mayer-Vietoris, the canonical map $\hocolim(\sP^n/\sH_{\sU},1)\xrightarrow{\epsilon}(\P^n/H,1)$
induced by the cover $\sU$ is an isomorphism in $\sH_\bullet(k)$.
The collection of quotient maps $U_{J^c}\to U_{J^c}/H\cap U_{J^c}$ or identity maps give the map $\gamma:\sP^n_{\sU,1}/\sH_{\sU_1}\to (\sP^n/\sH_{\sU},1)$.

We also have the functor $\sigma^n\G_m^{\wedge n}$. Identifying $U_{0\ldots n}$ with $\G_m^n$ via the coordinates $(x_1^0,\ldots, x_n^0)$, the quotient map $U_{0\ldots n}\cong \G_m^n\to \G_m^{\wedge n}$ extends canonically to the quotient map
$\delta:\sP^n_{\sU,1}/\sH_{\sU_1}\to \sigma^n\G_m^{\wedge n}$.
From our discussion on the isomorphism $\alpha$, we have the commutative diagram of isomorphisms in $\sH_\bullet(k)$
\begin{equation}\label{eqn:Comm2}
\xymatrix{
&(\P^n/H,1)\ar[drr]^\alpha\\
 \hocolim(\sP^n/\sH_{\sU},1)\ar[ur]^\epsilon&\hocolim\sP^n_{\sU,1}/\sH_{\sU_1}\ar[l]^-\gamma\ar[u]_{\epsilon_1} \ar[r]_-\delta&\hocolim\sigma^n\G_m^{\wedge n}\ar[r]_-\vartheta&\Sigma^n_n\G_m^{\wedge n}.
}
\end{equation}

Note that, for each $J\neq\0,[n]$, we have $\A_J\subset U_J$, 
since for $j\in J$, the intersection $\Pi_J\cap (X_j=0)$ is equal to $\Pi'_J$. Also, the map $\tilde\phi_J:\tilde{\sG}_n(J)\to \P^n$ has image contained in $U_{J^c}$. We define the map of functors
\[
\psi^T_*:\sG^T_n\to T_+\wedge (\sP^n/\sH_{\sU},1)
\]
as follows:  for $J\neq\0,[n]$, we use the map 
\[
(p_1,\phi^{T\prime}_J):\sG^{T\prime}_n(J)\to T\times U_{J^c}/(U_{J^c}\cap H)
\]
on $\sG^{T\prime}_n(J)$,  and the map
\[
T\times \A_{J^c} \xrightarrow{\id\times i_J} T\times U_{J^c} 
\]
induced by the inclusion $i_J:\A_{J^c}\hookrightarrow U_{J^c}$. One checks that these descend to a well defined map on the quotient
\[
\psi^T_J:\sG^T_n(J)\to T_+\wedge (\sP^n/\sH_{\sU},1)(J).
\]
For $J=\0$, we use 
\[
(\id_T,\phi^{T\prime}_\0):T\to  T\times U_{0\ldots n}/H\cap U_{0\ldots n}
\]
This gives us the commutative diagram of functors 
\[
\xymatrixcolsep{40pt}
\xymatrix{
\sG^T_n\ar[r]^-{\psi^T_*}\ar[dd]& T_+\wedge (\sP^n/\sH_{\sU},1)\\
&T_+\wedge \sP^n_{\sU,1}/\sH_{\sU_1}\ar[u]_{\id\wedge\gamma}\ar[d]^{\id\wedge \delta}\\
\sigma^n\sT\ar[r]_-{\id_T\wedge\sigma^n(\rho\circ f)}& T_+\wedge\sigma^n\G_m^{\wedge n}
}
\]
which induces the commutative diagram (in $\sH_\bullet(k)$) on the homotopy colimits
\[
\xymatrix{
\hocolim\sG^T_n\ar[r]^-{\Psi^T}\ar[d]_{\beta^T}& T_+\wedge \hocolim(\sP^n/\sH_{\sU},1)\ar[d]^{\id\wedge\vartheta\circ \delta\circ\gamma^{-1}}\\
\Sigma_s^n T_+\ar[r]_-{\Sigma_s^n(\id\wedge\rho\circ f)}& T_+\wedge\Sigma_s^n\G_m^{\wedge n}.
}
\]
Combining this with our diagram \eqref{eqn:Comm2} and noting that $\Phi^T=(\id\wedge\epsilon)\circ\Psi^T$ yields the commutative diagram
in $\sH_\bullet(k)$
\[
\xymatrix{
\hocolim\sG^T_n\ar[r]^-{\Psi^T}\ar@/^20pt/[rr]^-{\Phi^T}\ar[d]_{\beta^T}& T_+\wedge \hocolim(\sP^n/\sH_{\sU},1) \ar[d]^{\id\wedge(\vartheta\circ \delta\circ\gamma^{-1})}\ar[r]^-{\id\wedge\epsilon}& T_+\wedge(\P^n/H,1)\ar[dl]^-{\id\wedge\alpha}\\
\Sigma_s^nT_+\ar[r]_-{\Sigma_s^n(\id\wedge\rho\circ f)}& T_+\wedge\Sigma_s^n\G_m^{\wedge n},
}
\]
completing the proof.
\end{proof}

\section{Computing the collapse map}
We retain the notation from \S\S\ref{sec:PTCollapse}, \ref{sec:Susp1} and \ref{sec:Susp2}. Our task in this section  is to use lemma~\ref{lem:HauptLem} to give an explicit computation of $Q_F(w)$
as the $n$th suspension of a map $\rho_w:\Spec F_+\to w_+\wedge\G_m^{\wedge n}$, 
at least for $w$ an $F$-point of $\Delta^n\setminus\partial\Delta^n$. In general, we will need to take a further $\P^1$-suspension before desuspending, which we do in the next section.

For $F$ a finitely generated  field extension of $k$ and $w$ a closed point of $\Delta^n_F\setminus\partial \Delta^n_F$, we have the Pontryagin-Thom collapse map (definition~\ref{Def:PTCollapseMap})
\[
PT_F(w):\Sigma^n_s\Spec F_+\to (\P^n_{F(w)}/H_{F(w)},1).
\]
We have as well the map \eqref{eqn:CollapseMap} 
\[
Q_F(w):\Sigma^n_s\Spec F_+\to w_+\wedge\Sigma^n_s\G_m^{\wedge n}=w_+\wedge S^{2n,n}
\]
It follows from the definition of $MV_w$ \eqref{eqn:MVDef}, $PT_F(w)$ and $mv_w$ (definition~\ref{Def:PurityIsomv_w}) that
\begin{equation}\label{eqn:MainPlayer}
Q_F(w)=(\id_{w_+}\wedge\alpha)\circ PT_F(w),
\end{equation}
where we identify $ (\P^n_{F(w)}/H_{F(w)},1)$ with $w_+\wedge  (\P^n/H,1)$ and where $\alpha: (\P^n/H,1)\to \Sigma^n_s\G_m^{\wedge n}$ is the isomorphism \eqref{eqn:SuspIso1}. 

Consider an $F$-point $w:\Spec F\to \Delta^n$ of $\Delta^n$. Given elements $z_1,\ldots, z_n$ of $F^\times$, we have the corresponding map
\[
[z_1]\wedge_F\ldots\wedge_F[z_n]: \Spec F_+\to \Spec F_+\wedge \G_m^{\wedge n}
\]
given as the composition
\[
\Spec F_+\xrightarrow{\id\wedge(z_1,\ldots, z_n)}\Spec F_+\wedge(\G_m^n,1)\to\Spec F_+\wedge \G_m^{\wedge n}.
\]
We use the notation $\wedge_F$ to denote  the smash product for points $F$-schemes $(X,x)$, $(Y,y)$:
\[
(X,x)\wedge_F(Y,y):=X\times_FY/(X\times_Fy\vee x\times_FY),
\]
and note that $[z_1]\wedge_F\ldots\wedge_F[z_n]$ really is the $\wedge_F$-product of the maps $[z_i]$.

\begin{prop} \label{prop:Computation1} Take $w=(w_0,\ldots, w_n)\in (\Delta^n\setminus\partial\Delta^n)(F)$. Then
\[
Q_F(w)=\Sigma^n_s[-w_1/w_0]\wedge_F\ldots\wedge_F[-w_n/w_0].
\]
\end{prop}

\begin{proof}  We have for each $V_n$-scheme $T\to V_n$  the functor \eqref{eqn:GTDiag}; applying this construction for the
 morphism  $w: \Spec F\to V_n$,  gives us the functor 
 \[
 \sG^w_n:\square^{n+1}_{*<[n]}\to \Spc_\bullet(k).
 \]
  We recall the subschemes $U'_i$, $i=0,\ldots,n$ and $H$ of $\A^1\times\Delta^n$ from \S\ref{sec:Susp2}.

We note that $U'_0=1\times\Delta^n$, $H\cap U'_0$ is the face $t_0=0$, and that $U'_0\cap U'_i$ is the face $t_i=0$, for $i=1,\ldots, n$. Thus, collapsing the  $U'_i$, $i=1,\ldots, n$, $H\cap U'_0$ and all the $\A_J$ to a point, and sending $U_0'$ to $\Delta^n$ by the projection map gives a well defined morphism in $\Spc_\bullet(k)$, 
\[
a:\hocolim\sG^w_n\to  \Spec F _+\wedge\Delta^n/\partial\Delta^n,
\]
which is an isomorphism in $\sH_\bullet(k)$. In addition, we have the commutative diagram of isomorphisms in $\sH_\bullet(k)$
\begin{equation}\label{eqn:SuspDiag}
\xymatrix{
\hocolim\sG^w_n\ar[r]^-a_-\sim\ar[d]_{\beta^w}^\sim&\Spec F_+\wedge\Delta^n/\partial\Delta^n\\
\Sigma^n_s\Spec F_+\ar[ur]_{\sigma^F}^\sim,
}
\end{equation}
where $\sigma^F$ is the isomorphism \eqref{eqn:sigma} and  $\beta^w$ is the isomorphism \eqref{eqn:beta}.

Let 
\[
\tilde{r}_w:\Spec F_+\wedge(\P^n/H,1)\to \P^n_F/(\P^n_F\setminus \{(1:0:\ldots:0)\})
\]
be the composition of the isomorphism $\Spec F_+\wedge(\P^n/H,1)\cong (\P^n_F/H_F,1)$ followed by the quotient map 
$r_w:(\P^n_F/H_F,1)\to  \P^n_F/(\P^n_F\setminus \{(1:0:\ldots:0)\})$. It follows directly from the definition of the map $\Phi^w$ \eqref{eqn:**}
and the map $\phi_w$ \eqref{eqn:phi} that the diagram
\[
\xymatrix{
\hocolim\sG^w_n\ar[rr]^{\Phi^w}\ar[d]_a&&\Spec F_+\wedge  (\P^n/H,1)\ar[d]^{\tilde{r}_w}\\
\Spec F_+\wedge \Delta^n/\partial\Delta^n\ar[r]_{c_w}&\Delta^n_F/\Delta^n_F\setminus\{w\}\ar[r]_-{\phi_w}&\P_F^n/(\P_F^n\setminus\{(1\:0\:\ldots\:0)\}).
 }
 \]
 commutes. Combining this with the diagram \eqref{eqn:SuspDiag} and our description \eqref{eqn:PT} of $PT_F(w)$ gives us the commutative diagram
 \[
 \xymatrix{
\Sigma^n_s\Spec F_+\ar[rd]^-{PT_F(w)}\\
\hocolim\sG^w_n\ar[r]_-{\Phi^w}\ar[u]^{\beta^w}_\sim& \Spec F_+\wedge (\P^n/H,1) 
}
\]
But by lemma~\ref{lem:HauptLem},
\[
(\Sigma^n_s[-w_1/w_0]\wedge_F\ldots\wedge_F[-w_n/w_0])\circ \beta^w=(\id_{\Spec F_+}\wedge\alpha)\circ\Phi^w;
\]
since   $\beta^w$ is an isomorphism, this gives us 
\[
\Sigma^n_s[-w_1/w_0]\wedge_F\ldots\wedge_F[-w_n/w_0]=(\id_{\Spec F_+}\wedge\alpha)\circ PT_F(w).
\]
Our formula \eqref{eqn:MainPlayer} for $Q_F(w)$ completes the proof.
\end{proof}

\section{Transfers and $\P^1$-suspension}\label{sec:Transfer}

We now consider the general case of a closed point $w\in V_{nF}\subset\Delta^n_F$.  

Consider the map 
\begin{align*}
&j:\Delta^n\to \P^n\\
&j(t_0,\ldots, t_n):=(1:t_1,\ldots:t_n);
\end{align*}
$j$ is an open immersion, identifying $\Delta^n$ with $U_0$ and $V_n$ with  $U_{0\ldots n}\setminus H\subset\P^n$.

We define the {\em transfer map}
\[
Tr_F(w):S^{2n,n}\wedge\Spec F_+\to S^{2n,n}\wedge\Spec F(w)_+
\]
associated to a closed point $w\in \A^n_F$, separable over $F$,  as the composition
\begin{multline*}
S^{2n,n}\wedge\Spec F_+\xymatrix{&\ar[l]_{\alpha_\infty\wedge\id}^\sim}\P^n_F/H_{\infty F} \xrightarrow{c_{j{w}}} \P^n_F/\P^n_F\setminus\{j(w)\})
\hbox{$\xymatrix{&\ar[l]_{\bar{p}\circ\bar{j}}^\sim}$} \Delta^n_{F(w)}/(\Delta^n_{F(w)}\setminus \{w\})\\
 \xrightarrow{mv^\infty_w} \P^n_{F(w)}/H_{\infty F(w)}
\xymatrix{\ar[r]^{\alpha_\infty\wedge\id}_\sim&} S^{2n,n}\wedge\Spec F(w)_+.
\end{multline*}
The map $\bar{j}$ is induced from $j$, $\bar{p}$ is induced from the projection $p:\Delta^n_{F(w)}\to \Delta^n_F$, and $w\in \Delta^n_{F(w)}$ is the canonical lifting of $w\in \Delta^n_F$ to $\Delta^n_{F(w)}=w\times_F\Delta^n_F$. The map $\bar{p}\circ\bar{j}$ is
an isomorphism by Nisnevich excision (which is where we use the separability of $w$ over $F$). The map $mv^\infty_w$ is the Morel-Voevodsky purity isomorphism, where we use the generators $(t_1-w_1,\ldots, t_n-w_n)$ for $m_w$, together with the isomorphism
\[
r_w: \P^n_{F(w)}/H_{\infty F(w)}\to \P^n_{F(w)}/(\P^n_{F(w)}\setminus (1\:0\:\ldots\:0))
\]
induced by the identity on $\P^n_{F(w)}$. The map $\alpha_\infty$ is the isomorphism \eqref{eqn:Iso2}.

\begin{lem}\label{lem:Identity} Suppose that $w$ is in $V_n(F)$. Then $Tr_F(w)=\id$.
\end{lem}

\begin{proof} Let $w_0=(1\:0\:\ldots\:0)\in U_0\subset \P^n(k)$, giving us the purity isomorphism
\[
mv^\infty_{w_0}:U_{0}/(U_{0}\setminus w_0)\to \P^n/H_\infty
\]
defined via the choice of generators $(x_1,\ldots, x_n)$ for $m_{w_0}$.  The morphism
\[
(x\:x_1\:\ldots\: x_n): \A^1\times U_0\setminus 0\times w_0\to \P^n
\]
extends to an $\A^1$-bundle
\[
\pi:=(x\:x_1\:\ldots\: x_n):\Bl_{0\times w_0}\A^1\times U_0\to \P^n
\]
Furthermore, the restriction of $\pi$ to $1\times U_0$ extends to the identity map $\P^n\to \P^n$. From this, it follows that morphism in $\sH_\bullet(k)$,
\[
Tr_F(w_0): S^{2n,n}\to S^{2n,n}
\]
is the identity. On the other hand, let $T_{w}:\P^n_F\to \P^n_F$ be the automorphism extending translation by $w$ on $U_0$. Then $T_{w}$ acts by the identity on $\P^n_F/H_{\infty F}$, as we can extend $T_{w}$ to the $\A^1$ family of automorphisms $t\mapsto T_{tw}$ connecting $T_w$ with $\id$. Furthermore, $T_{-w}^*(x_1,\ldots, x_n)=(x_1-w_1,\ldots, x_n-w_n)$. From this it follows that 
\[
Tr_F(w)=T_w\circ Tr_F(w_0)\circ T_{-w}=\id.
\]
\end{proof}

\begin{prop} \label{prop:Computation2} Let $w=(w_0,\ldots, w_n)$ be a closed point of $V_{nF}$, separable over $F$. Then the $S^{2n,n}$-suspension of $Q_F(w)$:
\[
\id_{S^{2n,n}}\wedge Q_F(w):S^{2n,n}\wedge\Spec F_+\wedge S^{n,0}\to S^{2n, n}\wedge
w_+\wedge S^{2n,n}
\]
is equal to the map
$\Sigma^n_s\left((\id_{S^{2n,n}}\wedge [-w_1/w_0]\wedge_{F(w)}\ldots\wedge_{F(w)}[-w_n/w_0])\circ  Tr_F(w)\right)$.
 
\end{prop}

\begin{proof} Write $*_F$ for $\Spec F$.  We have the commutative diagram
\[
\xymatrixcolsep{0pt}
\xymatrix{
S^{2n,n}\wedge *_{F+}\wedge S^{n,0}\ar[d]_{\alpha_\infty^{-1}\wedge\sigma_F}\ar[rr]^{Tr_F(w)\wedge\id}&&S^{2n,n}\wedge w_+\wedge S^{n,0}\ar[d]^{\id\wedge\sigma}\\
\P^n/H_\infty\wedge*_{F+}\wedge\Delta^n/\partial\Delta^n\ar[dd]_{\id\wedge c_w}\ar[rd]^-{c_{j(w)}\wedge\id}&&S^{2n,n}\wedge w_+\wedge\Delta^n/\partial\Delta^n\ar[ddd]^{\id\wedge c_w}\\
&
\hbox to50pt{\hss$\P^n/(\P^n\setminus j(w))\wedge \Delta^n/\partial\Delta^n$\hss}\ar[d]^{\id\wedge c_w}\ar[ru]^-{(\alpha_\infty\circ m^\infty_{j(w)})\wedge\id\quad}&\\
\P^n/H_\infty\wedge\Delta^n_F/(\Delta^n_F\setminus w)\ar[r]^-{c_{j(w)}\wedge\id}\ar[d]_{\id\wedge mv_w}&
\hbox to 70pt{$\P^n/(\P^n\setminus j(w))\wedge \Delta^n/(\Delta^n\setminus w)$\hss}\ar[rd]_-{(\alpha_\infty\circ m^\infty_{j(w)})\wedge\id\quad\quad}\ar[dd]_{\id\wedge mv_w}\\
\hbox to100pt{$\P^n/H_\infty\wedge w_+\wedge(\P^n/H,1)$\hss}\ar[rd]_-{c_{j(w)}\wedge\id\quad}\ar[ddd]_{\alpha_\infty\wedge\id_{w_+}\wedge\alpha}
&&\hbox to75pt{\hss$S^{2n,n}\wedge w_+\wedge \Delta^n/(\Delta^n\setminus w)$}\ar[dd]^{\id\wedge mv_w}\\
&\hbox to130pt{\hss$\P^n/(\P^n\setminus j(w))\wedge w_+\wedge(\P^n/H,1)$\hss}\ar[rd]_-{(\alpha_\infty\circ m^\infty_{j(w)})\wedge\id\quad\quad}\\
&&S^{2n,n}\wedge w_+\wedge(\P^n/H,1)
\ar[d]^{\id\wedge\id_{w_+}\wedge\alpha}\\
S^{2n,n}\wedge w_+\wedge  S^{2n,n}\ar@{=}[rr]&&S^{2n,n}\wedge w_+\wedge  S^{2n,n};
}
\]
the commutativity follows either by definition of $Tr_F(w)$, or by identities of the form $(a\wedge 1)\circ(1\wedge b)=(1\wedge b)\circ (a\wedge 1)$, or (in the bottom pentagon) lemma~\ref{lem:Identity}.
The composition along the left-hand side is $\id_{S^{2n,n}}\wedge [(\id_{w_+}\wedge\alpha)\circ PT_F(w)]$; along the right-hand side we have $\id_{S^{2n,n}}\wedge   [(\id_{w_+}\wedge\alpha)\circ PT_{F(w)}(w)]$. Since $w$ is $F(w)$-rational, we may apply proposition~\ref{prop:Computation1} and our formula  \eqref{eqn:MainPlayer}  for $Q_F(w)$ to complete the proof.
\end{proof}

\section{Conclusion}\label{sec:Conclusion}
We can now put all the pieces together. For $E\in\Spt_{S^1}(k)$ fibrant, we have the associated fibrant object
$\Omega^n_TE:=\sHom_{\Spt(k)}(S^{2n,n}, E)$,
that is, $\Omega^n_TE$ is the presheaf
$(\Omega^n_TE)(X):= E(X_+\wedge S^{2n,n})$.
For each $n\ge1$, we have the canonical map
\[
\iota_n:E\to \Omega^n_T\Sigma^n_TE.
\]

Replacing $S^{2n,n}$ with $S^{n,n}=\G_m^{\wedge n}$, we have the fibrant object
\[
\Omega^n_{\G_m}E:=\sHom_{\Spt(k)}(S^{n,n}, E),
\]
defined as the presheaf $(\Omega^n_{\G_m}E)(X):= E(X_+\wedge \G_m^{\wedge n})$.

Given a closed point $w\in V_{nF}$, we define the map
\[
Tr_F(w)^*:\pi_m(\Omega^n_TE(w))\to \pi_m(\Omega^n_TE(F))
\]
as the composition
\begin{align*}
\pi_m(\Omega^n_TE(w))&=\Hom_{\SH_{S^1}(k)}(\Sigma^\infty_s(S^{2n,n}\wedge w_+),\Sigma^{-m}_sE)\\
&\xrightarrow{\Sigma^\infty_s(Tr_F(w)))^*}
\Hom_{\SH_{S^1}(k)}(\Sigma^\infty_s(S^{2n,n}\wedge \Spec F_+),\Sigma^{-m}_sE)\\
&=\pi_m(\Omega^n_TE(F)).
\end{align*}

\begin{Def} Take $E\in \SH_{S^1}(k)$ and let $n\ge1$ be an integer. An {\em $n$-fold $T$-delooping} of $E$ is an an object $\omega_T^{-n}E$ of $\SH_{S^1}(k)$ and an isomorphism $\iota_n:E\to \Omega^n_T\omega^{-n}_TE$ in $\SH_{S^1}(k)$.
\end{Def}

Given an  $n$-fold $T$-delooping  of $E$,  $\iota_n:E\to \Omega^n_T\omega^{-n}_TE$,    the map $Tr_F(w)^*$ for $\Omega^n_T\omega^{-n}_TE$ induces the ``transfer map"
\[
\iota_n^{-1}\circ Tr_F(w)^*\circ \iota_n: \pi_m(E(w))\to \pi_m(E(F)),
\]
which we write simply as $Tr_F(w)^*$.  

\begin{rems}\label{rems:Transfer} 1. The transfer map $Tr_F(w)^*:\pi_m(E(w))\to \pi_m(E(F))$ may possibly depend on the choice of $n$-fold $T$-delooping, we do not have an example, however.\\
\\
2. An $n-b$-fold $T$-delooping of $E$ gives rise to an $n$-fold $T$-delooping of $\Omega^b_{\G_m}E$. Thus, via the adjunction isomorphism
\[
\Pi_{a,b}E\cong \pi_a\Omega^b_{\G_m}E
\]
we have a transfer map 
\[
Tr_F(w)^*:\Pi_{a,b}E(w)\to \Pi_{a,b}E(F)
\]
for $w$ a closed point of $V_{nF}$, separable over $F$.\\
\\
3. If $E=\Omega_T^\infty\sE$ for some $\sE\in\SH(k)$, then $E$ admits canonical $n$-fold $T$-deloopings, namely
\[
\omega^{-n}_TE:=\Omega_T^\infty\Sigma_T^n\sE.
\]
Indeed, in $\SH(k)$, $\Sigma_T$ is the inverse to $\Omega_T$ and $\Omega_T^\infty$ commutes with $\Omega_T$.
\end{rems}

For a morphism $\phi:\Sigma^\infty_sw_+\to E$, we have the suspension $\Sigma^n_T\phi:\Sigma^n_T\Sigma^\infty_sw_+\to \Sigma^n_TE$, the composition
\[
\Sigma^n_T\phi\circ \Sigma^\infty_s Tr_F(w)^*:\Sigma^n_T\Sigma^\infty_s\Spec F_+\to \Sigma^n_TE
\]
and the adjoint morphism 
\[
(\Sigma^n_T\phi\circ \Sigma^\infty_s Tr_F(w)^*)':\Sigma^\infty_s\Spec F_+\to \Omega_T^n\Sigma^n_TE.
\]

Suppose we have an $n$-fold de-looping of $E$,  $\iota_n:E\to \Omega^n_T\omega^{-n}_TE$. This gives us the adjoint
\[
\iota'_n:\Sigma^n_TE\to \omega^{-n}_TE
\]
and
\[
\Omega^n_T\iota'_n:\Omega^n_T\Sigma^n_TE\to \Omega_T^n\omega^{-n}_TE.
\]
Let $\delta_n:E\to \Omega^n_T\Sigma^n_TE$ be the unit for the adjunction.

\begin{lem}\label{lem:Adjunction} 1. $\iota_n=\Omega^n_T\iota'_n\circ\delta_n$\\
\\
2. $\iota_n^{-1}\circ \Omega^n_T\iota'_n\circ (\Sigma^n_T\phi\circ \Sigma^\infty_s Tr_F(w))'=Tr_F(w)^*(\phi)$.
\end{lem}

\begin{proof} The two assertions follow from the universal property of adjunction.
\end{proof}

Before proving our main results, we show that the transfer maps respect the Postnikov filtration $F^*_\Tate\pi_mE$.

\begin{lem} \label{lem:TransferFilt} Suppose $E$ admits an $n$-fold $T$-delooping $\iota_n:E\to \Omega^n_T\omega^{-n}_TE$. Then for each finitely generated  field $F$ over $k$ and each closed point $w\in\A^n_F$ separable over $F$, we have
\[
Tr_F(w)^*(F^q_\Tate\pi_mE(w))\subset F^q_\Tate\pi_mE(F).
\]
\end{lem}

\begin{proof} Take $q\ge0$, and let $\tau_q:f_qE\to E$ be the canonical morphism. As above, let $\iota'_n:\Sigma^n_TE\to \omega^{-n}_TE$ be the adjoint of $\iota_n$ and let $\delta_n:E\to \Omega_T^n\Sigma_T^nE$ be the unit of the adjunction.  By lemma~\ref{lem:Adjunction}, we have the factorization of $\iota_n$ as
\[
E\xrightarrow{\delta_n}\Omega_T^n\Sigma_T^nE\xrightarrow{\Omega^n_T\iota'_n}\Omega_T^n\omega_T^{-n}E.
\]
This gives us the commutative diagram
\[
\xymatrix{
f_qE\ar[d]_{\delta_n}\ar[r]^{\tau_q}&E\ar[d]^{\iota_n}\\
\Omega^n_T\Sigma^n_Tf_qE\ar[r]_-{\tau'_q}& \Omega^n_T\omega^{-n}_TE,
}
\]
where $\tau_q':=\Omega^n_T\iota'_n\circ \Omega^n_T\Sigma^n_T\tau_q$.
Since $\iota_n:E\to \Omega^n_T\omega^{-n}_TE$ is an isomorphism, the composition 
\[
\iota_n\circ\tau_q:f_qE\to \Omega^n_T\omega^{-n}_TE
\]
satisfies the universal property of $f_q\Omega^n_T\omega^{-n}_TE\to \Omega^n_T\omega^{-n}_TE$.
By \cite[theorem 7.4.1]{LevineHC}, $\Omega^n_T\Sigma^n_Tf_qE$ is in $\Sigma^q_T\SH_{S^1}(k)$, hence there is a canonical morphism
\[
\theta:\Omega^n_T\Sigma^n_Tf_qE\to f_qE
\]
extending our first diagram to the commutative diagram
\[
\xymatrix{
f_qE\ar@<-3pt>[d]_{\iota_n}\ar[r]^{\tau_q}&E\ar@<3pt>[d]^{\iota_n}\\
\Omega^n_T\Sigma^n_Tf_qE\ar[r]_-{\tau'_q}\ar@<-3pt>[u]_\theta&\ar@<3pt>[u]^{\iota_n^{-1}} \Omega^n_T\omega^{-n}_TE.\\
}
\]
Using the universal property of $\tau_q$, we see that $\theta\circ\iota_n=\id_{f_qE}$, i.e., 
\[
\Omega^n_T\Sigma^n_Tf_qE=f_qE\oplus R
\]
and the restriction of $\tau'_q$ to $R$ is the zero map. We define the transfer map
\[
Tr_F(w)^*:\pi_mf_qE(w)\to \pi_mf_qE(F)
\]
by using the transfer map for $\Omega^n_T\Sigma^n_Tf_qE$ and this splitting.

The second diagram thus gives rise to the commutative diagram
\[
\xymatrix{
\pi_mf_qE(w)\ar[d]_{Tr_F(w)^*}\ar[r]^-{\tau_q}&\pi_mE(w)\ar[d]^{Tr_F(w)^*}\\
\pi_mf_qE(F)\ar[r]_-{\tau_q}&\pi_mE(F),
}
\]
which yields the result.
\end{proof}

\begin{rem} \label{rem:TransferFilt} One can define transfer maps in a more general setting, that is, for a closed point $w\in \A^n_F$ and any choice of parameters for $m_w\subset \sO_{\A^n,w}$. The same proof as used for lemma~\ref{lem:TransferFilt} shows that these more general transfer maps respect the filtration $F^*_\Tate\pi_mE$. \end{rem}

\begin{thm}\label{thm:Main1} Let $E\in\Spt(k)$ be fibrant, and let $F$ be a  field extensions of $k$.  \\
\\
1. For each $w=(w_0,\ldots, w_n)\in V_n(F)$, and each $\rho\in \pi_{0}\Omega^n_{\G_m}E(F)$, the element
\[
\rho\circ\Sigma^\infty_s([-w_1/w_0]\wedge_F\ldots\wedge_F[-w_n/w_0]):\Sigma^\infty_s\Spec F_+\to E
\]
is in $F^n_\Tate\pi_0E(F)$.\\
\\
2. Suppose that $E$ admits an $n$-fold $T$-delooping $\iota_n:E\to\Omega^n_T\omega^{-n}_TE$.  Then  for $w=(w_0,\ldots, w_n)$ a closed point of $V_{nF}$, separable over $F$, and $\rho_w\in \pi_{0}\Omega^n_{\G_m}E(w)$
\begin{equation}\label{eqn:Gen}
Tr_F(w)^*[\rho_w\circ\Sigma^\infty_s([-w_1/w_0]\wedge_F\ldots\wedge_F[-w_n/w_0])] 
\end{equation}
is in $F^n_\Tate\pi_0E(F)$.\\
\\
3. Suppose that $E$ admits an $n$-fold $T$-delooping $\iota_n:E\to\Omega^n_T\omega^{-n}_TE$. and  that $\Pi_{a,*}E=0$ for all $a<0$. Suppose further that $F$ is perfect. Then  $F^n_\Tate\pi_0E(F)$ is generated by elements of the form \eqref{eqn:Gen}, as $w$ runs over 
closed point of $V_{nF}$ and  $\rho_w$ over elements of $\pi_{0}\Omega^n_{\G_m}E(w)$.
\end{thm}

\begin{proof}  (1) follows directly from proposition~\ref{prop:Factorization} and proposition~\ref{prop:Computation1}, noting that
the isomorphism $\Omega_{\G_m}^nE\cong \Sigma^n_s\Omega^n_TE$ gives the identification 
\[
\pi_{-n} \Omega^n_TE(w)\cong \pi_0\Omega_{\G_m}^nE(w)\cong \Hom_{\SH_{S^1}(k)}(\Sigma^\infty_s w_+\wedge\G_m^{\wedge n},E).
\]

For (2), the fact that this element is in $F^n_\Tate\pi_0(E(F))$ follows from (1) and lemma~\ref{lem:TransferFilt}. 

For (3), that is, to see that these elements generate, take one of the generators  $\gamma:=\xi_w\circ\Sigma^\infty_sQ_F(w)$ of  
$F^n_\Tate \pi_0E(F)$, as given by proposition~\ref{prop:Factorization}, that is,   $w$ is a closed point of $V_{nF}$ and 
$\xi_w$ is in  $\pi_{-n}(\Omega_T^nE(w))=\pi_0(\Omega_{\G_m}^nE(w))$.  Since $F$ is perfect, $w$ is separable over $F$. Take the $n$-fold $T$-suspension of $\gamma$
\[
\Sigma^n_T\gamma:\Sigma^\infty_s(\Sigma_T^n \Spec F_+)\to \Sigma_T^nE,
\]
giving by adjunction and composition with $\Omega^n_T(\iota'_n)$ the morphism
\[
\Omega^n_T(\iota'_n)\circ (\Sigma^n_T\gamma)':\Sigma^\infty_s\Spec F_+\to \Omega_T^n\omega^{-n}E.
\]
It follows from the universal properties of adjunction that
\[
(\Sigma^n_T\gamma)'=\delta_n\circ\gamma,
\]
hence by lemma~\ref{lem:Adjunction} we have
\begin{equation}\label{eqn:***}
\Omega^n_T(\iota'_n)\circ (\Sigma^n_T\gamma)'=\Omega^n_T(\iota'_n)\circ\delta_n\circ\gamma=\iota_n\circ\gamma.
\end{equation}

Write
\[
\Sigma^n_T\gamma=(\Sigma^n_T\xi_w)\circ (\Sigma^\infty_s \Sigma^n_T Q_F(w)).
\]
By  proposition~\ref{prop:Computation2} we have
\[
\Sigma^n_T Q_F(w)=\Sigma^n_s\left(\Sigma^n_T[-w_1/w_0]\wedge_F\ldots\wedge_F[-w_n/w_0]\circ Tr_F(w)\right),
\]
and thus
\[
\Sigma^n_T\gamma=\Sigma^n_T(\xi_w\circ \Sigma^n_s[-w_1/w_0]\wedge_F\ldots\wedge_F[-w_n/w_0])\circ \Sigma^n_sTr_F(w).
\]

Using \eqref{eqn:***} and lemma~\ref{lem:Adjunction}, we have 
\begin{align*}
\iota_n\circ\gamma&=\Omega^n_T(\iota'_n)\circ (\Sigma^n_T\gamma)'\\
&=\Omega^n_T(\iota'_n)\circ [\Sigma^n_T(\xi_w\circ \Sigma^n_s[-w_1/w_0]\wedge_F\ldots\wedge_F[-w_n/w_0])\circ \Sigma^n_sTr_F(w)]'\\
&=\iota_n\circ Tr_F(w)^*(\xi_w\circ \Sigma^n_s[-w_1/w_0]\wedge_F\ldots\wedge_F[-w_n/w_0]),
\end{align*}
or
\[
 \gamma=Tr_F(w)^*[\rho_w\circ\Sigma^\infty_s([-w_1/w_0]\wedge_F\ldots\wedge_F[-w_n/w_0])].
\]
\end{proof}

We now assume that $E=\Omega^\infty_T\sE$ for some fibrant $T$-spectrum $\sE\in \Spt_T(k)$. Let $\mS_k$ denote the motivic sphere spectrum in $\Spt_T(k)$, that is, $\mS_k$ is a fibrant model of the suspension spectrum $\Sigma_T^\infty S^0_k$. We proceed to re-interpret theorem~\ref{thm:Main1} in terms of the canonical action of $\pi_0\Omega^\infty_T\mS_k(F)$ on $\pi_0E(F)$, which we now recall, along with some of the fundamental computations of Morel relating the Grothendieck-Witt group with endomorphisms of the motivic sphere spectrum.

We recall the Milnor-Witt sheaves of Morel, $\underline{K}^{MW}_n$ (see \cite[section 2]{MorelA1} for details). The graded sheaf $\underline{K}^{MW}_*:=\oplus_{n\in\Z}\underline{K}^{MW}_n$ has structure of a Nisnevich sheaf of associative graded rings. For a finitely generated  field $F$ over $k$,  the graded ring $K^{MW}_*(F):=\underline{K}^{MW}_*(F)$ has generators $[u]$ in degree 1, for $u\in F^\times$, and an additional  generator $\eta$ in degree $-1$, with relations
\begin{itemize}
\item $\eta[u]=[u]\eta$
\item $[u][1-u]=0$ (Steinberg relation)
\item $[uv]=[u]+[v]+\eta[u][v]$
\item $\eta(2+\eta[-1])=0$.
\end{itemize}

For later use, we note the following result:

\begin{lem}\label{lem:ExtendedSteinberg} Let $F$ be a field,  $u_1,\ldots, u_n\in F^\times$ with $\sum_iu_i=1$. Then $[u_1]\cdot\ldots\cdot [u_n]=0$ in $K^{MW}_0(F)$.
\end{lem}

\begin{proof} We use a number of relations in $K^{MW}_*(F)$, proved in \cite[lemma 2.5, 2.7]{MorelA1}. For $u\in F^\times$ we let $\<u\>$ denote the element $1+\eta[u]\in K^{MW}_0(F)$. We have the following relations, for $a, b\in F^\times$,
\begin{enumerate}
\item[i)] $K_0^{MW}(F)\text{ is central in }K^{MW}_*(F)$
\item[ii)] $[a][1-a]=0$\text{ for }$a\neq 1$ 
\item[iii)] $[ab]=[a]+\<a\>[b]$
\item[iv)] $[a^{-1}]=-\<a^{-1}\>[a]$
\item[v)] $[a][-a]=0$
\item[vi)] $[1]=0$.
\end{enumerate}
These yield the additional relation
\begin{enumerate}
\item[vii)] $[a][-a^{-1}]=0$.
\end{enumerate}
This follows by noting that
\begin{align*}
[a][-a^{-1}]&=[a](-\<-a^{-1}\>[-a])&\text{(iv)}\\
&=(-\<-a^{-1}\>)[a][-a]&\text{(i)}\\
&=0&\text{(v)}
\end{align*}

We prove the lemma by induction on $n$, the case $n=1$ being the relation (vi), the case $n=2$ the Steinberg relation (ii). Induction  reduces to showing
\[
[u][v]=[u+v][-v/u]\text{ for }u+v\neq0
\]
(in case  $u+v=0$ we use  (v)  to continue the induction).
For this, we have
\begin{align*}
[u][v]&=[u][v]+\<v\>[u][-u^{-1}]&\text{(vii)}\\
&=[u][v]+[u]\<v\>[-u^{-1}]&\text{(i)}\\
&=[u][-v/u]&\text{(iii)}\\
&=[u][-v/u]+\<u\>[1+v/u][-v/u]&\text{(ii)}\\
&=[u+v][-v/u]&\text{(iii)}
\end{align*}
\end{proof}

For $u\in F^\times$,  let $\<u\>$ denote the quadratic form $uy^2$ in the Grothendieck-Witt group $\GW(F)$. Sending $[u]\eta$ to $\<u\>-1$ extends to an isomorphism  \cite[lemma 2.10]{MorelA1}
\[
\vartheta_0:K^{MW}_0(F)\to \GW(F).
\]
In addition, for $n\ge1$, the image of $\times\eta^n:K^{MW}_n(F)\to K^{MW}_0(F)$
is an ideal $\eta^nK^{MW}_n(F)$ in $K^{MW}_0(F)$ and $\vartheta_0$ maps $\eta^nK^{MW}_n(F)$ isomorphically onto the ideal $I(F)^n$, where $I(F)\subset \GW(F)$ is the augmentation ideal of quadratic forms of virtual rank zero.

For each $u\in F^\times$, we have the corresponding morphism 
\[
[u]:\Spec F_+\to \G_m
\]
We have as well the canonical projection $\eta':\A^2\setminus\{0\}\to\P^1$. Using a construction similar to the one we used to show that $\P^2/H\cong \Sigma^2_s\G_m^{\wedge 2}$, one constructs a canonical isomorphism in $\sH_\bullet(k)$, $(\A^2\setminus\{0\},1)\cong \Sigma^1_s\G_m^{\wedge 2}$,
and thus $\eta'$ yields the morphism
\[
\eta:\Sigma^1_s\G_m^{\wedge 2}\to \Sigma^1_s\G_m
\]
in $\sH_\bullet(k)$. 

For $E, F\in\Spt_{S^1}(k)$, let $\underline{\Hom}(E,F)$ denote the Nisnevich sheaf associated to the presheaf
\[
U\mapsto \Hom_{\SH_{S^1}(k)}(U_+\wedge E,F).
\]

We have the fundamental theorem of Morel:

\begin{thm}[\hbox{\cite[corollary 3.43]{MorelA1}}] \label{thm:MorelMain} Suppose $\Char k\neq 2$. Let $m,p,q\ge0$, $n\ge2$ be integers. Then sending $[u]\in K^{MW}_1(F)$ to the morphism $[u]$ and sending $\eta\in K^{MW}_{-1}(F)$ to the morphism $\eta$ yields isomorphisms
\[
\Hom_{\sH_\bullet(k)}(\Spec F_+\wedge S^m\wedge\G_m^{\wedge p}, S^n\wedge\G_m^{\wedge q})\cong 
\begin{cases}0&\text{ if }m<n\\K^{MW}_{q-p}(F)&\text{ if }m=n\text{ and } q>0.\end{cases}
\]
\end{thm}
As we will be relying on Morel's theorem, we assume for the rest of the paper that the characteristic of $k$ is different from two.

Passing to the $S^1$-stabilization, theorem~\ref{thm:MorelMain} gives
\begin{align}\label{align:MorelIso}
&\Pi_{0,p}\Sigma^\infty_s\G_m^{\wedge q}=\underline{K}^{MW}_{q-p}&\text{ for } p\ge0, q\ge1, \\
&\Pi_{a,p}\Sigma^\infty_s\G_m^{\wedge q}=0&\text{ for } p\ge0, q\ge1, a<0.\notag
\end{align}

Passing to the $T$-stable setting, Morel's theorem gives
\begin{align}\label{eqn:MorelIso}
&\pi_{p,p}\Sigma_{\G_m}^q\mS_k\cong \underline{K}^{MW}_{q-p}&\text{ for } p,q\in\Z  \\
&\pi_{a+p,p}\Sigma_{\G_m}^q\mS_k=0&\text{ for } p,q\in\Z, a<0.\notag
\end{align}

Composition of morphisms gives us the (right) action of the bi-graded sheaf of rings $\pi_{*,*}\mS_k$ on $\pi_{*,*}\sE$ for each $T$-spectrum $\sE$, and thus, the action of $\underline{K}^{MW}_{-*}$ on $\pi_{*,*}\sE$. If we let $E$ be the $S^1$-spectrum $\Omega_T^\infty\sE$, then
$\Pi_{a,b}E=\pi_{a+b,b}\sE$
for all $b\ge0$. Thus,  via lemma~\ref{lem:DegreeShift}(2) we thus have the right multiplication
\[
\Pi_{a,b-m}E\otimes \underline{K}^{MW}_{-m}\to \Pi_{a,b}E.
\]

Let $\sI\subset  \underline{K}^{MW}_0$ be the sheaf of augmentation ideals. The $\underline{K}^{MW}_{-*}$-module structure on $\Pi_{a,*}E$ gives us the filtration $F^{MW}_n\Pi_{a,b}E$ of $\Pi_{a,b}E$, defined by
\[
F^{MW}_n\Pi_{a,b}E:=\text{im}[\Pi_{a,n}E\otimes \underline{K}^{MW}_{n-b}\to \Pi_{a,b}E];\quad n\ge 0.
\]

\begin{lem}\label{lem:FIltProps}  Suppose $E=\Omega_T^\infty\sE$ for some $\sE\in\SH(k)$. For  integers $n, b, p\ge0$,  with $n-p, b-p\ge0$, the adjunction isomorphism $\Pi_{a,b}E\cong \Pi_{a,b-p}\Omega^p_{\G_m}E$ induces an isomorphism
\[
F^{MW}_{n}\Pi_{a,b}E\cong F^{MW}_{n-p}\Pi_{a,b-p}\Omega^p_{\G_m}E.
\]
\end{lem}

\begin{proof} This follows easily from the fact that the adjunction isomorphism 
 \[
 \Pi_{a,*}E\cong \Pi_{a,*-p}\Omega^p_{\G_m}E
 \]
 is a  $\underline{K}^{MW}_*$-module isomorphism.
\end{proof}

\begin{Def} Let $E=\Omega_T^\infty\sE$ for some $\sE\in \SH(k)$, $F$ a field extension of $k$. Take integers $a, b, n$ with $n, b\ge0$. Following remark~\ref{rems:Transfer}(2), we have the transfer maps 
\[
Tr_F(w):\Pi_{a,b}E(F(w))\to \Pi_{a,b}E(F)
\]
for each closed point $w\in V_{nF}$, separable over $F$. \\
\\
1. Let  $F_{n}^{MW\widehat{\ }_{Tr}}\Pi_{a,b}E(F)$ denote the subgroup of $\Pi_{a,b}E(F)$ generated by elements of the form
\[
Tr_F(w)^*(x);\quad x\in F^{MW}_{n}\Pi_{a,b}E(F(w))
\]
as $w$ runs over closed points of $V_{nF}$, separable over $F$.\\
\\
2. Let $[\Pi_{a,b}E\cdot\sI^n]^{\widehat{\ }_{Tr}}(F)$ denote the subgroup of $\Pi_{a,b}E(F)$ generated by elements of the form
\[
Tr_F(w)^*(x\cdot y);\quad x\in \Pi_{a,b}E(F(w)),  y\in I(F(w))^n, 
\]
as $w$ runs over closed points of $V_{nF}$, separable over $F$.
\end{Def}

\begin{rem}\label{rem:ProjForm}  It follows directly from the definitions that, for $w$ a closed point of $V_{nF}$, $x\in K^{MW}_{n-b}(F)$, $y\in \Pi_{a,n}E(F(w))$, we have
\[
Tr_F(w)^*(y\cdot p^*x)=Tr_F(w)^*(y)\cdot x,
\]
where $p^*x\in K^{MW}_{n-b}(F(w))$ is the extension of scalars of of $x$. In particular,  $[\Pi_{a,b}E\cdot\sI^n]^{\widehat{\ }_{Tr}}(F)$ is a $K^{MW}_0(F)$-submodule of   $\Pi_{a,b}E(F)$ containing $\Pi_{a,b}E(F)I(F)^n$.
\end{rem}

\begin{thm}\label{thm:Main2} Let $k$ be a perfect field of characteristic $\neq2$.  Let $E=\Omega_T^\infty\sE$ for some $\sE\in \SH(k)$ with $\Pi_{a,b}\sE=0$ for all $a<0$, $b\ge0$. Let $n> p\ge0$ be integers and let $F$ be a perfect field extension of $k$.   Then 
\[
F^n_\Tate\Pi_{0,p}E(F)=F_{n}^{MW\widehat{\ }_{Tr}}\Pi_{0,p}E(F).
\]
For $p\ge n\ge0$, we have the identity of sheaves $F^n_\Tate\Pi_{0,p}E= \Pi_{0,p}E$.
\end{thm}

\begin{proof} First suppose $n> p$.  By lemma~\ref{lem:DegreeShift}   and lemma~\ref{lem:FIltProps}, we reduce to the case $p=0$. 

The fact that we have an inclusion of  $K^{MW}_0(F)$-submodules of $\Pi_{0,0}E(F)$,
\[
F^n_\Tate\Pi_{0,0}E(F)\subset F_{n}^{MW\widehat{\ }_{Tr}}\Pi_{0,0}E(F),
\]
follows from theorem~\ref{thm:Main1}. Indeed, as $F$ is perfect,  each element of the form \eqref{eqn:Gen} is of the form $Tr_F(w)(\rho_w\cdot z)$, with
$\rho_w\in \Pi_{0,n}E(w)$, $z\in K^{MW}_n(F(w))$, hence in $F_{n}^{MW\widehat{\ }_{Tr}}\Pi_{0,0}E(F)$.

To show the other inclusion, it suffices by lemma~\ref{lem:TransferFilt} and  theorem~\ref{thm:Main1} to show that, for each field $K$ finitely generated over $k$, the elements $[-u_1/u_0]\cdot\ldots\cdot[-u_n/u_0]$, with $(u_0,\ldots, u_n)\in V_n(K)$, generate  $K^{MW}_n(K)$ as a module over $K^{MW}_0(K)$.

We note that the map sending $(u_0,\ldots, u_n)$ to $(1/u_0, -u_1/u_0,\ldots,-u_n/u_0)$ is an involution of $V_n$, so it suffices to show that the elements $[u_1]\cdot\ldots\cdot[u_n]$, with $(u_0,\ldots, u_n)\in V_n(K)$, generate. 

Sending $(u_0,\ldots, u_n)$ to $(u_1,\ldots, u_n)$ identifies $V_n$ with $(\A^1\setminus\{0\})^n\setminus H$. But by definition $K^{MW}_n(K)$ is generated by elements 
$[u_1]\cdot\ldots\cdot[u_n]$ with $u_i\in K^\times$; it thus suffices to show that $[u_1]\cdot\ldots\cdot[u_n]=0$ in $K^{MW}_n(K)$ if $\sum_iu_i=1$; this is lemma~\ref{lem:ExtendedSteinberg}.

If $p\ge n\ge0$, the universal property of $f_nE\to E$ gives us the isomorphism for $U\in \Sm/k$ 
\[
 \Hom_{\SH_{S^1}(k)}(\Sigma_s^\infty \Sigma_{\G_m}^pU_+, E)\cong  \Hom_{\SH_{S^1}(k)}(\Sigma_s^\infty \Sigma_{\G_m}^pU_+, f_nE),
\]
since $\Sigma_s^\infty \Sigma_{\G_m}^pU_+$ is in $\Sigma^p_T\SH_{S^1}(k)$ for $U\in \Sm/k$. As these groups of morphisms define the presheaves whose respective sheaves are  $\Pi_{0,p}E(F)$ and $\Pi_{0,p}f_nE$, the map $\Pi_{0,p}f_nE\to \Pi_{0,p}E$ is an isomorphism,  hence $F^n_\Tate\Pi_{0,p}E= \Pi_{0,p}E$.
\end{proof}
 
\begin{rem}\label{rem:ExtendedTransfer} The reader may object that the collection of transfer maps used to define $F_{n}^{MW\widehat{\ }_{Tr}}\Pi_{0,p}E(F)$ is rather artificial. However, the fact that the general transfer maps mentioned in remark~\ref{rem:TransferFilt} respect the filtration $F^*_\Tate\pi_mE$, together with theorem~\ref{thm:Main2}, shows that, if we were to allow arbitrary transfer maps in our definition of $F_{n}^{MW\widehat{\ }_{Tr}}\Pi_{0,p}E(F)$, we would arrive at the same subgroup of $\Pi_{0,m}E(F)$.
\end{rem}

Our main result for a $T$-spectrum, theorem~\ref{thm:main0}, follows easily from theorem~\ref{thm:Main2}:
\begin{proof}[Proof of \hbox{ theorem~\ref{thm:main0}}] Using lemma~\ref{lem:DegreeShift}, we reduce to the case $p=0$. Essentially the same argument as used at the end of the proof of theorem~\ref{thm:Main2} proves the part of 
theorem~\ref{thm:main0} for $n\le 0$.

If $n>0$, then for $b\ge 0$,  we have
\begin{align*}
&\pi_{a,b}\sE\cong \pi_{a,b}\Omega^\infty_T\sE&\text{(lemma~\ref{lem:DegreeShift})}\\
&\pi_{a,b}f_n\sE\cong \pi_{a,b}\Omega^\infty_Tf_n\sE\cong  \pi_{a,b}f_n\Omega^\infty_T\sE&\text{(lemma~\ref{lem:DegreeShift}) and } \eqref{eqn:SliceIso}
\end{align*}
 Thus, in case $n>0$, theorem~\ref{thm:main0} for $\sE$ is equivalent to 
 theorem~\ref{thm:Main2} for $\Omega_T^\infty\sE$, completing the proof.
 \end{proof}

Finally, we can prove our main result for the motivic sphere spectrum, theorem~\ref{thm:main}.  Let $\sE=\Sigma^q_{\G_m}\mS_k$. Then  Morel's isomorphism \eqref{eqn:MorelIso} and lemma~\ref{lem:DegreeShift} give
\[
\Pi_{a,b}\Omega_T^\infty\sE=\begin{cases}\underline{K}^{MW}_{q-b}&\text{ for }a=0, b\ge0\\0&\text{ for }a<0, b\ge0.\end{cases}
\]

\begin{thm}\label{thm:Main3} Let $k$ be a perfect field of characteristic $\neq2$.\\
 1. For all $n>p\ge0$, $q\in\Z$,  and all perfect field extensions $F$ of $k$,  we have
\[
F^n_\Tate\Pi_{0,p}\Omega_T^\infty\Sigma^q_{\G_m}\mS_k(F)= K^{MW}_{q-p}(F)I(F)^N\subset K^{MW}_{q-p}(F), 
\]
where $N=N(n-p, n-q):=\max(0,\min(n-p, n-q))$.  In particular,
\[
F^n_\Tate\pi_{0,0}\mS_k(F)= I(F)^n\subset \GW(F).
\]
2. For $n\le p$, we have the identity of sheaves $F^n_\Tate\Pi_{0,p}\Omega_T^\infty\Sigma^q_{\G_m}\mS_k= \underline{K}^{MW}_{q-p}$.\\
\\
3. In case $k$ has characteristic zero, we have the identity of sheaves
\[
F^n_\Tate\Pi_{0,p}\Omega_T^\infty\Sigma^q_{\G_m}\mS_k=  \underline{K}^{MW}_{q-p}\sI^N\subset \underline{K}^{MW}_{q-p}. 
\]
with $N$ as above.
\end{thm}

\begin{proof} Let $N$ be as defined in the statement of the theorem. We first note (3) follows from (1), in fact, from (1)  for all fields extensions $F$ finitely generated over $k$.  Indeed, 
$F^n_\Tate\Pi_{0,p}\Omega_T^\infty\Sigma^q_{\G_m}\mS_k$  is the image of the map
\[
\Pi_{0,p}f_n\Omega_T^\infty\Sigma^q_{\G_m}\mS_k\to \Pi_{0,p}\Omega_T^\infty\Sigma^q_{\G_m}\mS_k
\]
induced by the canonical morphism $f_n\Omega_T^\infty\Sigma^q_{\G_m}\mS_k\to \Omega_T^\infty\Sigma^q_{\G_m}\mS_k$. By results of Morel \cite[theorem 3 and lemma 5]{MorelConn}, both homotopy sheaves are strictly $\A^1$-invariant sheaves of abelian groups. But the category of strictly $\A^1$-invariant sheaves of abelian groups is abelian \cite[lemma 6.2.13]{MorelConn}, hence $F^n_\Tate\Pi_{0,p}\Omega_T^\infty\Sigma^q_{\G_m}\mS_k$ is also strictly $\A^1$-invariant. It follows, e.g., from Morel's isomorphism
\[
\pi_0\Omega^\infty_T\Sigma^m_{\G_m}\mS\cong \pi_{-m,-m}\mS\cong \underline{K}^{MW}_m
\]
that the sheaves $ \underline{K}^{MW}_m$ are strictly $\A^1$-invariant; as $\underline{K}^{MW}_{q-p}\sI^N$ is the image of the map
\[
\times\eta^M\:\underline{K}^{MW}_{q-p+M}\to \underline{K}^{MW}_{q-p},
\]
where $M=N$ if $q-p\ge0$, $M=p-q+N$ if $q-p<0$,
it follows that  $\underline{K}^{MW}_{q-p}\sI^N$ is strictly $\A^1$-invariant as well. Our assertion follows from the fact  that a strictly $\A^1$-invariant sheaf $\sF$ is zero if and only $\sF(k(X))=0$ for all $X\in\Sm/k$, which in turn is an easy consequence of \cite[lemma 3.3.6]{MorelLec}.

Next, suppose $n-p\le 0$. Then $N=0$ and
\begin{align*}
F^n_\Tate\Pi_{0,p}\Omega_T^\infty\Sigma^q_{\G_m}\mS_k&=F^{n-p}_\Tate\Pi_{0,0}\Omega^p_{\G_m}\Omega_T^\infty\Sigma^q_{\G_m}\mS_k
&\text{(lemma~\ref{lem:DegreeShift})} \\
&=\Pi_{0,0}\Omega^p_{\G_m}\Omega_T^\infty\Sigma^q_{\G_m}\mS_k&(n-p<0)\\
&=\Pi_{0,p}\Omega_T^\infty\Sigma^q_{\G_m}\mS_k&\text{(adjunction)}\\
&=\underline{K}^{MW}_{q-p}&\text{(Morel's theorem)}
\end{align*}
proving (2); we  may thus assume $n-p>0$.

By \eqref{eqn:MorelIso},  we may apply theorem~\ref{thm:Main2}, which tells us that $F^n_\Tate\Pi_{0,p}\Omega_T^\infty\Sigma^q_{\G_m}\mS_k(F)$ is the  subgroup of $\Pi_{0,p}\Omega_T^\infty\Sigma^q_{\G_m}\mS_k(F)=K^{MW}_{q-p}(F)$ generated by elements 
 of the form  $Tr_F(w)^*(y\cdot x)$ with 
 \begin{align*}
 &y\in \Pi_{0,n}\Omega_T^\infty\Sigma^q_{\G_m}\mS_k(F(w))=K^{MW}_{q-n}(F(w))\\
 &x \in K^{MW}_{n-p}(F(w)).
 \end{align*}
 
 Suppose that $n-q< 0$, so $N=0$. Then $q-n\ge 0$ and $n-p>0$, and thus the product map
 \[
\mu_{n-p,q-n}: K^{MW}_{n-p}(F(w))\otimes K^{MW}_{q-n}(F(w))\to K^{MW}_{q-p}(F(w))=\Pi_{0,p}\Omega_T^\infty\Sigma^q_{\G_m}\mS_k(F(w))
 \]
 is surjective. Since the map $Tr_F(w)$ is an isomorphism for $w\in V_n(F)$, we see that
 \[
 F^n_\Tate\Pi_{0,p}\Omega_T^\infty\Sigma^q_{\G_m}\mS_k(F)=\Pi_{0,p}\Omega_T^\infty\Sigma^q_{\G_m}\mS_k(F).
 \]
 
 Suppose $n-q\ge 0$. Then
  \[
 \times\eta^{n-q}:K^{MW}_0(F(w))\to  K^{MW}_{q-n}(F(w))
 \]
 is surjective.   If $n-p\ge n-q$, then the image of   $\mu_{n-p,q-n}$ is the same as the image of the triple product
 \[
 K^{MW}_{q-p}(F(w))\otimes  K^{MW}_{n-q}(F(w))\otimes K^{MW}_{q-n}(F(w))\to K^{MW}_{q-p}(F(w));
 \]
 as the image of
 \[
 \mu_{n-q,q-n}:K^{MW}_{n-q}(F(w))\otimes K^{MW}_{q-n}(F(w))\to  K^{MW}_{0}(F(w))
 \]
 is $I(F(w))^{n-q}$, we see that the image of $\mu_{n-p,q-n}$ is $K^{MW}_{q-p}(F(w))I(F(w))^{n-q}$ and thus
 \[
F_\Tate^n\Pi_{0,p}\Omega_T^\infty\Sigma^q_{\G_m}\mS_k(F)=
[\Pi_{0,p}\Omega_T^\infty\Sigma^q_{\G_m}\mS_k\sI^N]^{\widehat{\ }_{Tr}}(F).
\]
Similarly, if $n-q\ge n-p$, then the image of   $\mu_{n-p,q-n}$ is the same as the image of the triple product
 \[
 K^{MW}_{q-p}(F(w))\otimes  K^{MW}_{n-p}(F(w))\otimes K^{MW}_{p-n}(F(w))\to K^{MW}_{q-p}(F(w))
 \]
 which is $K^{MW}_{q-p}(F(w))I(F(w))^{n-p}$. Thus
 \[
F_\Tate^n\Pi_{0,p}\Omega_T^\infty\Sigma^q_{\G_m}\mS_k(F)=
[\Pi_{0,p}\Omega_T^\infty\Sigma^q_{\G_m}\mS_k\sI^N]^{\widehat{\ }_{Tr}}(F)
\]
in this case as well.

Thus, to complete the proof, it suffices to show that, for $w$ a closed point of $V_{nF}$, and $N\ge0$ an integer, we have
\begin{equation}\label{eqn:ToShow}
Tr_F(w)^*\left(K^{MW}_{q-p}(F(w))I(F(w))^N\right)\subset {K}^{MW}_{q-p}(F)I(F)^N.
\end{equation}
First suppose that $q-p\ge0$. Take a closed point $w\in V_{nF}$ and elements $x_1,\ldots, x_N\in F(w)^\times$, $y\in {K}^{MW}_{q-p}(F(w))$. We have
\begin{align*}
Tr_F(w)^*(y\cdot [x_1]\eta\cdot\ldots\cdot [x_N]\eta)&=Tr_F(w)^*(y\cdot [x_1] \cdot\ldots\cdot [x_N]\eta^N)\\
&=Tr_F(w)^*(y\cdot [x_1] \cdot\ldots\cdot [x_N])\cdot \eta^N.
\end{align*}
where we  use remark~\ref{rem:ProjForm} in the last line. Since $q-p\ge0$,  ${K}^{MW}_{q-p}(F)I(F)^N$ is the image in ${K}^{MW}_{q-p}(F)$ of the map
\[
-\times\eta^N:{K}^{MW}_{q-p+N}(F)\to K^{MW}_{q-p}(F),
\]
which verifies \eqref{eqn:ToShow}.

In case $q-p<0$, write $y=y_0\eta^{p-q}$, with $y_0\in K^{MW}_0(F(w))$. As above, we have
\[
Tr_F(w)^*(y\cdot [x_1]\eta\cdot\ldots\cdot [x_N]\eta)=Tr_F(w)^*(y_0\cdot [x_1]\cdot\ldots\cdot [x_N])\cdot \eta^{p-q+N},
\]
which is in $\eta^{p-q}\cdot [K^{MNW}_N(F)\eta^N]=K^{MW}_{q-p}(F)I(F)^N$, as desired.
\end{proof}

Theorem~\ref{thm:Main3} yields the main result for the $S^1$-spectra $\Sigma^\infty_s\G_m^{\wedge q}$ by using  the $S^1$-stable consequences of Morel's unstable computations, theorem~\ref{thm:MorelMain}.

\begin{cor} \label{cor:S1Stable}  Let $k$ be a perfect field of characteristic $\neq2$. \\
1.  For all $n>p\ge0$, $q\ge1$, and all perfect field extensions $F$ of $k$, we have
\[
F^n_\Tate\Pi_{0,p}\Sigma^\infty_s\G_m^{\wedge q}(F)=K^{MW}_{q-p}(F)I(F)^{N(n-p,n-q)}\subset K^{MW}_{q-p}(F),
\]
with $N(n-p,n-q)$ as in theorem~\ref{thm:Main3}. \\
\\
2. For $n\le p$, we have $F^n_\Tate\Pi_{0,p}\Sigma^\infty_s\G_m^{\wedge q}= \Pi_{0,p}\Sigma^\infty_s\G_m^{\wedge q}$.\\
\\
3. If $\Char k=0$, we have the identity of sheaves
\[
F^n_\Tate\Pi_{0,p}\Sigma^\infty_s\G_m^{\wedge q}=\underline{K}^{MW}_{q-p}\sI^{N(n-p,n-q)}\subset \underline{K}^{MW}_{q-p}.
\]
\end{cor}

\begin{proof} As in the proof of theorem~\ref{thm:Main3}, it suffices to prove (1).

The main point is that Morel's unstable computations show that the $\G_m$-stabilization  map
\begin{multline*}
\Hom_{\SH_{S^1}(k)}(\Sigma^m_s\Sigma^\infty_s\G_m^{\wedge p}\wedge\Spec F_+, \Sigma^\infty_s\G_m^{\wedge q})\\\to
\Hom_{\SH_{S^1}(k)}(\Sigma^m_s\Sigma^\infty_s\G_m^{\wedge p+1}\wedge\Spec F_+, \Sigma^\infty_s\G_m^{\wedge q+1})
\end{multline*}
is an isomorphism for all $m\le0$, $p\ge0$ and $q\ge1$.

Let $E(p,q)=\Omega^p_{\G_m}\Sigma^\infty_s\G_m^{\wedge q}$, and let 
\[
E(q-p)=\Omega^\infty_T\Sigma^{-p}_{\G_m}\Sigma^\infty_T\G_m^{\wedge q}=\Omega^\infty_T\Sigma^{q-p}_{\G_m}\mS_k.
\]
Then 
\[
\pi_aE(p,q)=\Pi_{a,p}\Sigma^\infty_s\G_m^{\wedge q}.
\]
Thus $\Pi_{a,*}E(p,q)=0$ for $m<0$ and so we may apply proposition~\ref{prop:Factorization} to give generators of the form $\xi_w\circ\Sigma^\infty_sQ_F(w)$ for 
\[
F^{n-p}_\Tate\Pi_{0,0}\Omega^p_{\G_m}\Sigma^\infty_s\G_m^{\wedge q}(F)
=F^n_\Tate\Pi_{0,p}\Sigma^\infty_s\G_m^{\wedge q}(F).
\]
But $\xi_w$ is in 
\[
\pi_{-n+p}\Omega_T^{n-p}E(p,q)(w)=\pi_{0,n-p} E(p,q)(w).
\]
Similarly, we have generators $\xi'_w\circ\Sigma^\infty_sQ_F(w)$ for $F^{n-p}_\Tate\pi_0E(p-q)(F)$, with 
\[
\xi'_w\in  \pi_{0,n-p} E(p-q)(w).
\]
But the stabilization map
\[
\pi_{0,n-p} E(p,q)(w)\to \pi_{0,n-p}E(p+1,q+1)(w)
\]
is an isomorphism, and hence we have an isomorphism from the generators for $F^{n-p}_\Tate\pi_0E(p,q)(F)$ to the generators for 
\[
F^{n-p}_\Tate\pi_0E(q-p)(F)=\colim_mF^{n-p}_\Tate\pi_0E(p+m,q+m)(F).
\]
As the map
\[
\pi_0E(p,q)(F)\to \pi_0E(q-p)(F)=K^{MW}_{q-p}(F)
\]
is an isomorphism, it follows that the surjection
\[
F^{n-p}_\Tate\pi_0E(q-p)(F)\to F^{n-p}_\Tate\pi_0E(q-p).
\]
is an isomorphism as well. By theorem~\ref{thm:Main3}, we have
\[
F^{n-p}_\Tate\pi_0E(q-p)=K^{MW}_{q-p}(F)I(F)^N\subset K^{MW}_{q-p}(F),
\]
 completing the proof.
\end{proof}

Theorem~\ref{thm:Main3} also gives us the $T$-stable version
\begin{cor}\label{cor:MainTStable}  Let $k$ be a perfect field of characteristic $\neq2$. For  $n,p,q\in\Z$,  and $F$ a perfect field extensions of $k$, we have
\[
F^n_\Tate\pi_{p,p}\Sigma^q_{\G_m}\mS_k(F)= K^{MW}_{q-p}(F)I(F)^{N(n-p,n-q)}\subset K^{MW}_{q-p}(F)
\] 
For $n\le p$, we have $F^n_\Tate\pi_{p,p}\Sigma^q_{\G_m}\mS_k=  \underline{K}^{MW}_{q-p}$.  If $\Char k=0$, we have
\[
F^n_\Tate\pi_{p,p}\Sigma^q_{\G_m}\mS_k=  \underline{K}^{MW}_{q-p}\sI^{N(n-p,n-q)}\subset \underline{K}^{MW}_{q-p}.
\]
\end{cor}

\begin{proof} Using  lemma~\ref{lem:DegreeShift}   and lemma~\ref{lem:FIltProps} as in the proof of  theorem~\ref{thm:Main2} we have
\[
F^n_\Tate\pi_{p,p}\Sigma^q_{\G_m}\mS_k=F^{n-p+r}_\Tate\pi_{r,r}\Sigma^{q-p+r}_{\G_m}\mS_k
\]
for all integers $r$. As our assertion is also stable under this shift operation,  we may assume that $p,q\ge0$. We note that $\mS_k$ is in $\SH^\eff(k)$, hence so are all $\Sigma^q_{\G_m}\mS_k$ for $q\ge0$, and thus 
\[
F^n_\Tate\pi_{p,p}\Sigma^q_{\G_m}\mS_k=\pi_{p,p}\Sigma^q_{\G_m}\mS_k
\]
for $n<0$, $p,q\ge0$. The truncation functors  $f_n$, $n\ge0$,  on $\SH(k)$ and $\SH_{S^1}(k)$ commute with $\Omega_T^\infty$, and $\pi_{a,p}\Omega_T^\infty\sE=\pi_{a,p}\sE$ for $\sE\in \SH(k)$, $p\ge0$. This reduces us to computing  computing $F^n_\Tate\pi_{p,p}\Omega_T^\infty\Sigma^q_{\G_m}\mS_k$ for $n,p, q\ge0$, which is theorem~\ref{thm:Main3}.
\end{proof}

\section{Epilog: Convergence questions}\label{sec:Convergence} Voevodsky has stated a conjecture \cite[conjecture 13]{VoevOpen} that would imply that for  $\sE=\Sigma^\infty_TX_+$, $X\in\Sm/k$, the Tate Postnikov tower is convergent in the following sense: for all $a,b,n\in\Z$, one has
\[
\cap_mF^m_\Tate\pi_{a,b}f_n\sE=0. 
\]
Our computation of $F^n_\Tate\pi_{p,p}\Sigma_T^\infty\G_m^{\wedge q}$ gives some evidence for this convergence conjecture.

\begin{prop}\label{prop:converge}  Let $k$ be a perfect field with $\Char k\neq2$.  For all $p,q\ge0$, and all perfect field extensions $F$ of $k$, we have
\[
\cap_nF^n_\Tate\pi_{p,p}\Sigma^\infty_T\G_m^{\wedge q}(F)=0.
\]
\end{prop}

\begin{proof} In light of theorem~\ref{thm:Main3},  the assertion is that the $I(F)$-adic filtration on $K^{MW}_{q-p}(F)$ is separated. By \cite[th\'eor\`eme 5.3]{MorelWitt}, for $m\ge0$, $K^{MW}_m(F)$ fits into a cartesian square of $\GW(F)$-modules
\[
\xymatrix{
K^{MW}_m(F)\ar[r]\ar[d]&K^M_m(F)\ar[d]^{Pf}\\
I(F)^m\ar[r]_-q&I(F)^m/I(F)^{m+1},
}
\]
where $K^M_m(F)$ is the Milnor $K$-group, $q$ is the quotient map and $Pf$ is the map sending a symbol $\{u_1,\ldots,u_m\}$ to the class of the Pfister form $\<\<u_1,\ldots,u_m\>\>$ mod $I(F)^{m+1}$.  For $m<0$, $K^{MW}_m(F)$ is isomorphic to the Witt group of $F$, $W(F)$, that is, the quotient of $\GW(k)$ by the ideal generated by the hyperbolic form $x^2-y^2$. Also, the map $\GW(F)\to W(F)$ gives an isomorphism of $I(F)^r$ with its image in $W(F)$ for all $r\ge1$. Thus
\[
K^{MW}_m(F)I(F)^n=\begin{cases}I(F)^n\subset W(F)&\text{ for }m<0, n\ge0\\
I(F)^{n+m}\subset  \GW(F)&\text{ for }m\ge0, n\ge1.
\end{cases}
\]
The fact that $\cap_nI(F)^n=0$ in $W(F)$ or equivalently in $\GW(F)$ is a theorem of Arason and Pfister \cite{ArasonPfister}.
\end{proof}

\begin{rems} 1. The proof in  \cite{MorelWitt} that $K^{MW}_m(F)$ fits into a cartesian square as above relies the Milnor conjecture.\\
\\
2. Voevodsky's conjecture [{\it loc. cit.}] asserts the convergence for a wider class of objects in $\SH(k)$ than just the $T$-suspension spectra of smooth $k$-schemes. The selected class is the triangulated category generated by $\Sigma^n_T\Sigma^\infty_TX_+$, $X\in\Sm/k$, $n\in\Z$ and the taking of direct summands. However, as pointed out to me by Igor Kriz, the convergence fails for this larger class of objects. In fact, take $\sE$ to be the Moore spectrum $\mS_k/\ell$ for some prime $\ell\neq 2$. Since $\Pi_{a,q}\mS_k=0$ for $a<0$,   proposition~\ref{prop:Connected} shows that $\Pi_{a,q}f_n\mS_k=0$ for $a<0$, and  thus we have the right exact sequence for all $n\ge0$
\[
\pi_{0,0}f_n\mS_k\xrightarrow{\times \ell} \pi_{0,0}f_n\mS_k\to  \pi_{0,0}f_n\sE\to0.
\]
In particular, we have 
\[
F^n_\Tate\pi_{0,0}\sE(k)=im\left(F^n_\Tate\pi_{0,0}\mS_k(k)\to \pi_{0,0}\mS_k(k)/\ell\right)=im\left(I(k)^n\to \GW(k)/\ell\right).
\]
Take $k=\R$. Then $\GW(\R)=\Z\oplus\Z$, with virtual rank and virtual index giving the two factors. The augmentation ideal $I(\R)$ is thus isomorphic to $\Z$ via the index and it is not hard to see that $I(\R)^{n}=(2^{n-1})\subset\Z=I(\R)$. Thus $\pi_{0,0}\sE=\Z/\ell\oplus \Z/\ell$ and the filtration $F^n_\Tate\pi_{0,0}\sE$  is constant, equal to $\Z/\ell=I(\R)/\ell$, and is therefore not separated.

The convergence property is thus not a ``triangulated" one in general, and therefore seems to be quite subtle. However, if the $I$-adic filtration on $\GW(F)$ is finite (possibly of varying length depending on $F$) for all finitely generated $F$ over $k$, then our computations (at least in characteristic zero) show that the filtration $F^*_\Tate\pi_{p,p}\Sigma^\infty_T\G_m^{\wedge q}$ is at least locally finite, and thus has better triangulated properties; in particular, for $\ell\neq2$, 
\[
\pi_{0,0}(\mS_k/\ell)=\Z/\ell,\ F^n_\Tate\pi_{0,0}(\mS_k/\ell)=0\text{ for }n>0,
\]
as the augmentation ideal in $\GW(F)$ is purely two-primary torsion, and $\sI\pi_{0,0}\mS_k/\ell=0$. One can therefore ask if Voevodsky's convergence conjecture is true if one assumes the finiteness of the $I(F)$-adic filtration on $\GW(F)$ for all finitely generated fields $F$ over $k$. 
\end{rems}


\begin{thebibliography}{99}
\bibitem{ArasonPfister}
Arason, J. K. and Pfister, A., {\sl
Beweis des Krullschen Durchschnittsatzes f\"ur den Wittring}. Invent. Math. {\bf 12} (1971), 173--176. 

\bibitem{BousfieldKan}
 Bousfield, A. K. and Kan, D. M.{\bf  Homotopy limits, completions and localizations}. Lecture Notes in Mathematics, {\bf 304}. Springer-Verlag, Berlin-New York, 1972.

\bibitem{GoerssJardine}
Goerss, P. G.; Jardine, J. F., {\sl Localization theories for simplicial presheaves}. Canad. J. Math. {\bf 50} (1998), no. 5, 1048--1089. 

 \bibitem{Jardine}
 Jardine, J. F., {\sl Simplicial presheaves}, J. Pure Appl. Algebra {\bf 47} (1987), no. 1, 35--87

 \bibitem{Jardine2}
Jardine, J.F., {\sl Motivic symmetric spectra},
{\em Doc. Math.} 5 (2000), 445--553.

\bibitem{LevineHC}
Levine, M.,  {\sl The homotopy coniveau tower}. J Topology {\bf 1} (2008) 217--267.

\bibitem{LevineChow} Levine, M., {\sl Chow's moving lemma in $\A^1$-homotopy theory}.
$K$-theory {\bf 37} (1-2)  (2006) 129-209.


\bibitem{MorelA1}
Morel, F., {\sl $\A^1$-Algebraic topology over a field}, preprint 19.11.2006. http://www.mathematik.uni-muenchen.de/~morel/A1homotopy.pdf

\bibitem{MorelConn}
Morel, F., {\sl The stable $\A^1$-connectivity theorems}, $K$-theory, 2005, {\bf 35}, pp 1-68.

\bibitem{MorelWitt}
 Morel, F. {\sl Sur les puissances de l'id\'eal fondamental de l'anneau de Witt}. Comment. Math. Helv. {\bf 79} (2004), no. 4, 689--703. 


\bibitem{MorelLec} 
Morel, F.,  {\sl An introduction to $\mathbb A\sp 1$-homotopy theory}.  {\em Contemporary developments in algebraic $K$-theory} 357--441, ICTP Lect. Notes, {\bf XV}, Abdus Salam Int. Cent. Theoret. Phys., Trieste, 2004.


\bibitem{MorelVoev}  Morel, F. and Voevodsky, V., {\sl $\A^1$-homotopy theory of schemes},
Inst. Hautes \'Etudes Sci. Publ. Math. {\bf 90} (1999), 45--143.

\bibitem{Neeman}
Neeman, A,
{\bf Triangulated categories}.
Annals of Mathematics Studies, {\bf 148}. Princeton University Press, Princeton, NJ, 2001. 

\bibitem{VoevOpen}
Voevodsky, V., {\sl Open problems in the motivic stable homotopy theory. I}.  {\em Motives, polylogarithms and Hodge theory, Part I (Irvine, CA, 1998)}, 3--34, Int. Press Lect. Ser., {\bf 3}, I, Int. Press, Somerville, MA, 2002. 

\bibitem{VoevSlice}  
Voevodsky, V. {\sl A possible new approach to the motivic spectral sequence for algebraic $K$-theory}, {\em  Recent progress in homotopy theory (Baltimore, MD, 2000)} 371--379, Contemp. Math., {\bf 293} (Amer. Math. Soc., Providence, RI, 2002). 
\end{thebibliography}
\end{document}